\documentclass{article}
\usepackage{latexsym}
\usepackage{amsfonts}
\usepackage{amsmath}
\usepackage{amsfonts}
\usepackage{epsfig}
\usepackage{stmaryrd}
\usepackage{amssymb}
\usepackage{pxfonts}
\usepackage{wasysym}
\usepackage[font={small}]{caption}
\newtheorem{theorem}{Theorem}

\usepackage{algorithmic}
\usepackage{algorithm}

 \author{ Francisco Bernal \footnotemark[1] \ \footnotemark[2]}

\title{Trust-Region Methods for Nonlinear Elliptic Equations with Radial Basis
Functions}
\date{\vspace{-5ex}}
\begin{document}
\maketitle

\renewcommand{\thefootnote}{\fnsymbol{footnote}}

\footnotetext[1]{INESC-ID$\backslash$IST, TU Lisbon. Rua Alves Redol 9, 1000-029 Lisbon, Portugal.}
\footnotetext[2]{Center for Mathematics and its Applications,
Department of Mathematics, Instituto Superior T\'ecnico. 
Av. Rovisco Pais 1049-001 Lisbon, Portugal. ({\tt francisco.bernal@ist.utl.pt})}
\renewcommand{\thefootnote}{\arabic{footnote}}

%\begin{flushleft}

%--------------------------------------
\begin{abstract}
%--------------------------------------
We consider the numerical solution of nonlinear elliptic boundary value
problems with Kansa's method. We derive analytic formulas for the Jacobian and
Hessian of the resulting nonlinear collocation system and exploit them within
the framework of the trust-region algorithm. This ansatz is tested on semilinear, quasilinear and fully nonlinear elliptic PDEs (including Plateau's problem, Hele-Shaw flow and the Monge-Amp\`ere equation) with excellent
results. The new approach distinctly outperforms previous ones based on linearization or finite-difference Jacobians.  
\end{abstract}


%--------------------------------------------------------------------------
\section{Introduction}
\label{S:Introduction}
%--------------------------------------------------------------------------

\subsection{RBF interpolation}\label{SS:Interpolation}

Given the scalar data ${u_1,\ldots,u_N}$ on a set (called {\em pointset}) of distinct points ${\bf x}_1,\dots,{\bf x}_N\in{\mathbb R}^d$ (called {\em centres}), the RBF interpolant is defined as
\begin{equation}
\label{F:RBF_Interpolant}
{\tilde u}({\bf x})= \sum_{j=1}^N \alpha_j \phi(\parallel {\bf x}-{\bf x}_j\parallel),
\end{equation} 
%where the coefficients $\alpha_1,\ldots,\alpha_N$ are to be determined and
where the function $\phi(r):[0,\infty)\shortrightarrow{\mathbb R}$ is the chosen radial basis function (RBF). A few popular RBFs are shown in Table \ref{T:RBFs}. Thoughout this paper, $||\cdot||$ is always the 2-norm. The coefficients $\alpha_1,\ldots,\alpha_N$ are determined by collocation 

\begin{equation}
\label{F:Interpolation_system}
\left[
\begin{array}{ccc}
\phi(||{\bf x}_1-{\bf x}_1||) &	\ldots 	& \phi(||{\bf x}_1-{\bf x}_N||) \\
\vdots & \ddots & \vdots \\
\phi(||{\bf x}_N-{\bf x}_1||) &	\ldots 	& \phi(||{\bf x}_N-{\bf x}_N||)
\end{array}
\right]
\left[
\begin{array}{c}
\alpha_1 \\ \vdots \\ \alpha_N
\end{array}
\right]=
\left[
\begin{array}{c}
u_1 \\ \vdots \\ u_N
\end{array}
\right],
\end{equation}
or, more compactly, $[\phi]{\vec\alpha}={\vec u}$. Thanks to the radial argument of $\phi$, $[\phi]$--called the {\em RBF interpolation matrix}--is symmetric. Guaranteed non-singularity of $[\phi]$ depends on the RBF $\phi$ being strictly conditionally positive definite (SCPD)--i.e. bound to yield positive-definite $[\phi]$. For instance, in Table \ref{T:RBFs} all the RBFs are SCPD except for the multiquadric, where the RBF interpolant needs to be augmented with a constant to yield a positive definite $[\phi]$ \cite{Fasshauer_libro,Wendland_libro}. 

\begin{table}[h!]
\begin{footnotesize}
\[\begin{array}{lclll}
\multicolumn{5}{c}{\textrm{RBFs used in this paper}}\\
\multicolumn{5}{c}{ }\\
\hline\noalign{\smallskip}
RBF &  \phi(r) &  \textrm{notation} & \textrm{support} & \textrm{convergence rate}\\
\noalign{\smallskip}\hline\noalign{\smallskip}
\textrm{multiquadric} & \sqrt{r^2+c^2} & \textrm{MQ($c$)} & r\leq \infty & \textrm{spectral} \\
\textrm{inverse multiquadric} & 1/\sqrt{r^2+c^2} & \textrm{IMQ($c$)} & r\leq \infty & \textrm{spectral} \\
%\textrm{Mat\'ern} & C(r/c)^{(\alpha-d)/2}K_{(\alpha-d)/2}(r/c) & \textrm{MATERN($\alpha,c$)} & r\leq\infty &\textrm{spectral} \\
\textrm{Mat\'ern} & (r/c)^{(\alpha-d)/2}K_{(\alpha-d)/2}(r/c) & \textrm{MATERN($\alpha,c$)} & r\leq\infty & \textrm{spectral} \\
%\textrm{Wendland $C^4$ (in ${\mathbb R}^2$)} & [1-(r/L)]_+^{s+2}P(r/L,s) & \textrm{WC4($L$)} & r\leq L & \textrm{algebraic} \\
\textrm{Wendland $C^4$} & [1-(r/L)]_+^{s+2}P(r/L,s) & \textrm{WC4($L$)} & r\leq L & \textrm{algebraic} \\
\hline
\end{array}\]
\end{footnotesize}
%\caption{$d$ is the space dimension (${\bf x}\in{\mathbb R}^d$). In MATERN($\alpha,c$), the constant $C$ is $C=c^d\pi^{d/2}2^{(d+\alpha-2)/2}\Gamma(\alpha/2)$, where $K_{\nu}(t)$ is the modified Bessel function of the second kind. In WC4($L$), $[f(t)]_+=0$ if $t\geq 1$, $s=3+\lfloor d/2\rfloor$, and $P(t,s)=(s^2+4s+3)t^2+(3s+6)t+3$.}
\caption{$d$ is the space dimension (${\bf x}\in{\mathbb R}^d$). In MATERN($\alpha,c$), $K_{\nu}(t)$ is the modified Bessel function of the second kind. In WC4($L$), $[f(r)]_+=0$ if $r\geq L$, $s=3+\lfloor d/2\rfloor$, and $P(t,s)=(s^2+4s+3)t^2+(3s+6)t+3$ (WC4 works up to $d=3$).}
\label{T:RBFs}
\end{table}

\subsection{Kansa's method}\label{SS:Kansa}
In 1990, Kansa adapted this approach to the solution of linear boundary value problems (BVPs) \cite{Kansa90a,Kansa90b}. Consider the elliptic BVP
\begin{equation}
\label{F:Linear_PDE}
\left\{
\begin{array}{lr}
{\cal L}^{PDE}u({\bf x})= f, & \textrm{ if ${\bf x}\in\Omega$} \\
{\cal L}^{BC}u({\bf x})= g, & \textrm{ if ${\bf x}\in\partial\Omega$},
\end{array}
\right.
\end{equation}
where $\Omega$ is a bounded domain in ${\mathbb R}^d$, $d\geq 1$, $u:\Omega\rightarrow{\mathbb R}$ is smooth, and ${\cal L}^{PDE}$ and ${\cal L}^{BC}$ are the interior and boundary linear operators, respectively. % For example, for Poisson's equation with source term $f$ and Dirichlet data $g$, ${\cal L}=\nabla^2$ and ${\cal B}=I$ (the identity operator). 
Kansa's idea was to discretize $\Omega\cup\partial\Omega$ into a pointset $\Xi_N=\{{\bf x}_i\}_{i=1}^N$, and look for an approximation ${\tilde u}$ to $u$  with an RBF interpolant like (\ref{F:RBF_Interpolant}). Without loss of generality, we can assume that the first $M$ nodes in $\Xi_N$ belong to the interior of $\Omega$ and the last $N-M$ are discretizing its boundary. By linearity, collocation of (\ref{F:Linear_PDE}) on that interpolant leads to
%\begin{eqnarray}
%\label{F:RBF_Linear_System}
%{\cal L}{\tilde u}({\bf x}_i)=L\sum_{j=1}^N \alpha_j\phi_{ij}=\sum_{j=1}^N\alpha_jL\phi_{ij}, \qquad 1\leq i \leq M \nonumber\\
%{\cal B}{\tilde u}({\bf x}_i)=B\sum_{j=1}^N \alpha_j\phi_{ij}=\sum_{j=1}^N\alpha_jB\phi_{ij}, \qquad M< i \leq N,
%\end{eqnarray}  
%where $L\phi_{ij}=(L\phi)(\parallel {\bf x}_i-{\bf x}_j \parallel_2)$, etc. %Taking the MQ RBF for illustration, 
%\begin{equation}
%\nabla^2\sqrt{\parallel {\bf x}_i-{\bf x}_j\parallel_2^2 + c^2}=,
%\end{equation}
%but notice that $L\phi$ need not be radial if $L$ is not isotropic. For instance, $\frac{\partial\phi}{\partial x}({\bf x})= $, and even $\frac{\partial\phi}{\partial x}_{ij}=-\frac{\partial\phi}{\partial x}_{ji}$. In compact form, (\ref{F:RBF_Linear_System}) reads 
%The coefficients $\vec\alpha$ are now determined by
\begin{equation}
[{\cal L}\phi]{\vec\alpha}:=
\left[
\begin{array}{c}
[{\cal L}^{PDE}\phi]_{\Omega} \\
{[{\cal L}^{BC}\phi]}_{\partial\Omega}
\end{array}
\right]{\vec\alpha}:=
\left[
\begin{array}{cccccc}
{\cal L}^{PDE}\phi_{11} & \ldots & {\cal L}^{PDE}\phi_{1N} \\
\vdots & \ddots & \vdots \\ 
{\cal L}^{PDE}\phi_{M1}   & \ldots & {\cal L}^{PDE}\phi_{MN} \\
{\cal L}^{BC}\phi_{M+1,1} & \ldots & {\cal L}^{BC}\phi_{M+1,N} \\
\vdots & \ddots & \vdots \\
{\cal L}^{BC}\phi_{N1} & \ldots & {\cal L}^{BC}\phi_{NN} \\
\end{array}
\right]
{\vec\alpha}
=
\left[
\begin{array}{c}
f({\bf x}_1) \\ \vdots \\ f({\bf x}_M) \\ g({\bf x}_{M+1}) \\ \vdots \\ g({\bf x}_N)
\end{array}
\right].
\label{F:Kansa}
\end{equation}

%where ${\cal L}^{PDE}\phi_{ij}=({\cal L}\phi)(\parallel {\bf x}_i-{\bf x}_j \parallel)$, etc.
(Check Section \ref{SS:Notation} for the notation.) This method for solving PDEs has many appealing features: it is meshless,very easy to code, appropriate for high-dimensional PDEs (thanks to the radial argument of the RBFs, which is dimension-blind) and--as long as the solution is smooth--enjoys exponential convergence with respect to the fill distance of the pointset $\Xi_N$ (for many RBFs at least, see Table \ref{T:RBFs}). For a complete exposition, the reader is referred to \cite{Fasshauer_libro}. Regarding solvability, conditions which guarantee that the {\em differentiation matrix} in (\ref{F:Kansa}) be nonsingular have not yet been established. In fact, there are crafted examples which yield a singular matrix \cite{Hon2001}), but such cases should be exceedingly rare, as also confirmed by years of praxis. On the other hand, Kansa's method may lead to very ill-conditioned matrices, meaning that only pointsets with up to a few thousands of nodes can be used before the matrix in (\ref{F:Kansa}) becomes numerically singular. %Larger problems can still be tackled--at the expense of sacrificing the exponential accuracy--by using compactly supported RBFs such as WC4 in Table \ref{T:RBFs}. Finally, let us mention that a very promising new idea is the partition of unity-RBF method (RBF-PUM) for PDEs \cite{Larsson2013,Larsson2015}.   
Larger problems can be tackled by using compactly supported RBFs such as WC4 in Table \ref{T:RBFs} (at the expense of sacrificing spectral convergence), by the RBF-QR method \cite{RBF-QR} (for some RBFs), and/or by using the novel RBF-partition of unity method \cite{Larsson2015}.%\cite{Larsson2013,Larsson2015}.

\subsection{Nonlinear equations}\label{SS:Nonlinear} 

Extending Kansa's method to nonlinear equations is straightforward. %For definiteness, we will restrict ourselves in this paper to nonlinear elliptic BVPs, although most of the discussion could be carried over to other kind of nonlinear problems. % such as state-dependent delay differential equations \cite{Yodelay}. Let us specifically consider the nonlinear elliptic BVP compactly written as
Let us introduce the following compact notation for a nonlinear elliptic BVP:
\begin{equation}
\label{F:Nonlinear_Elliptic_BVP}
{\cal W}[{\bf x},u({\bf x}),Du({\bf x})]=0 \Rightarrow \left\{
\begin{array}{ll}
{\cal W}^{PDE}=0, & \textrm{ if ${x}\in\Omega$}\\
{\cal W}^{BC}=0, & \textrm{ if ${x}\in\partial\Omega$,}
\end{array}
\right.
\end{equation} 
%where $\Omega\subset {\mathbb R}^d$, $d\geq 1$, %${\cal W}$ represents either the BC or PDE operator depending on whether $\bf x\in\partial\Omega$ or $\bf x\in\Omega\backslash\partial\Omega$, respectively, and
where $Du({\bf x})$ is shorthand notation for any kind of derivatives present in (\ref{F:Nonlinear_Elliptic_BVP}), such as $\partial/\partial x,\nabla^2$, etc. Collocation of (\ref{F:RBF_Interpolant}) on (\ref{F:Nonlinear_Elliptic_BVP}) leads to the nonlinear system %of equations
\begin{equation}
\label{F:Nonlinear_System}
W_i({\vec\alpha}):= {\cal W}[{\bf x}_i,{\tilde u}({\bf x}_i),D{\tilde u}({\bf x}_i)] = 0, \qquad 1\leq i \leq N.
\end{equation} 
%(Note that some of the collocation equations may still be linear, for instance if Dirichlet BCs are present on a portion of the boundary). 
A root ${\vec\alpha}_*$ of (\ref{F:Nonlinear_System})--i.e. $\{W_i({\vec\alpha}_*)=0\}_{i=1}^N$ or simply ${\vec W}=0$%, check Section (\ref{SS:Notation})
--represents an RBF solution ${\tilde u}({\vec\alpha}_*)$ of the BVP (\ref{F:Nonlinear_Elliptic_BVP}). Even if the nonlinear BVP (\ref{F:Nonlinear_Elliptic_BVP}) has one unique solution, the meshless discretization (\ref{F:Nonlinear_System}) may have none, one, multiple or infinitely many roots, regardless of the fact that the system is square. Therefore, it is not evident that collocation is the best approach to RBF representations of solutions to nonlinear BVPs, especially given that least-squares RBF approximations have been found to be preferable to strict collocation in other contexts \cite{Ling2006,Platte2006}. Interestingly enough, we have found apparently unique strict roots in every well-conditioned square RBF collocation system arising from the various PDEs in our numerical experiments, provided that the domain discretization is reasonable enough. %In Section \ref{S:RBF_trust} we suggest a possible way to attack the questions of solvability and uniqueness of (\ref{F:Nonlinear_System}). 

%Anyways, our numerical experiments show that the existence of a strict root to (\ref{F:Nonlinear_System}) is rather unimportant, for two reasons. One is that, numerically, it may be hard to tell the RBF approximation associated to a root from that very closely approximating--but not interpolating--the PDE at the PDE collocation points, as long as the BCs are properly enforced. The second is that, in nonlinear equations, the error of the RBF approximation stops to drop shortly after the interpolation residual (i.e., the residual to the nonlinear PDE on non-collocation points) stagnates. Beyond that point, there is no real gain from iterating a solver of nonlinear systems forward to strict collocation (for such solvers are all iterative and require an initial guess ${\vec\alpha}_0$ (or ${\tilde u}_0$) to kick off the iterations).

\subsection{Rootfinding approach} 

We are seeking a root ${\vec\alpha}_*$ of (\ref{F:Nonlinear_System}). The most well-known rootfinding algorithm is Newton's method for systems \cite{Nocedal99}, which proceeds as:
\begin{equation}
\label{F:Newton_Method}
{\vec\alpha_{k+1}}= {\vec\alpha_{k}} - J_k^{-1}{\vec W}_k
\end{equation} 
where ${\vec W}_k=\Big(\,W_1(\vec\alpha_{k}),\ldots,W_N(\vec\alpha_{k})\,\big)^T$ and $J_k$ is the Jacobian evaluated at ${\vec \alpha}_k$. % (the entries of $J_k$ are $(J_{ij})_k=\partial W_i({\vec\alpha k})\backslash \partial \alpha_j$). 
%The step
%\begin{equation}
%\label{F:Full_Newton_step}
%\vec\gamma_N:={\vec\alpha_{k+1}} - {\vec\alpha_{k}}= - J_k^{-1}{\bf W}_k  
%\end{equation}
%is called the ${\textit full Newton step}$. 
Like all the solvers considered in this paper, Newton's method requires an initial guess ${\vec\alpha}_0$ (which may be the interpolation coefficients of a guess function  ${\tilde u}_0$) to kick off the iterations. One advantage of Newton's method is that, if ${\vec\alpha}_0$ is close enough to a root ${\vec\alpha}_*$, if $\det J({\vec\alpha}_*)\neq 0$, and if every Jacobian in the sequence (\ref{F:Newton_Method}) is well conditioned enough, then $\{{\vec \alpha}_k\}$ will converge to ${\vec\alpha}_*$ quadratically. However, the sequence may not be convergent at all if ${\vec\alpha}_0$ is not a good enough guess.
%Therefore, Newton's method is \textit{locally} but not \textit{globally} convergent, meaning that it is not guaranteed to find a root (if a root exists) from \textit{any} initial guess ${\vec\alpha}_0$. 
%Moreover, Newton's method breaks down if it comes accross a (numerically) singular Jacobian. Powell \cite{Powell1970} provides an influential example of a simple system where Newton's iterations are attracted to a singular Jacobian. The potentially erratic behaviour of Newton's method deems it inadequate. %Instead, most practical rootfinding algorithms are globally convergent modifications of it. Among them, 
In this paper, we analyze and advocate the trust-region algorithm (TRA) for nonlinear RBF collocation. %(And refer to the combination of both for solving nonlinear BVPs as the RBF/TRA method.) The TRA is considered superior in terms of robustness and speed of convergence \cite{Nocedal99}. Since there are many implementations of it %(like in Matlab's optimization toolbox)
%, it can be readily used to solve system (\ref{F:Nonlinear_System}). 
A TRA approach was implicitly used (via Matlab's \verb|fsolve|) to solve Navier-Stokes equations in \cite{Chinchapatnam2007}, with numerical Jacobians constructed via finite differences. Finite difference Jacobians are very expensive to construct and not as accurate as analytic ones; in particular for RBF collocation, they become numerically unstable long before. In this paper we derive--for the first time, to the best of our knowledge--analytic formulas for the Jacobian and Hessian of a wide range of nonlinear operators. Not only do they substantially improve the performance of the RBF/TRA method, but also enable a root of the nonlinear system to be found where the previous approaches fail, in the first place. Moreover, they also might offer theoretical insight into the root structure of the system.%\newline            

Let us briefly mention two directions that we have not pursued further. When all nonlinearities are made up of sums and products of derivatives--such as $u\nabla^2 u+(\partial u/\partial x)(\partial u/\partial y)$, for instance--RBF collocation gives rise to a system of polynomials in $\alpha_1,\ldots,\alpha_N$. In principle, the complete root structure of such a system could be revealed in the framework of Groebner bases \cite{Buchberger1985}, using for instance SINGULAR. (In practice, however, typical RBF systems are too large, and probably too ill-conditioned, to be tackled this way.) Another interesting method of tackling polynomial systems is homotopy/continuation.% \cite{Verschelde1999}.%\newline        

\subsection{Operator-Newton (linearization) approach} 
An alternative way of solving nonlinear elliptic BVPs with RBFs 
%yet without dealing with nonlinear systems of equations and rootfinding
is the operator-Newton method introduced by Fasshauer \cite{Fasshauer02}.  
The idea is to recast the original nonlinear BVP into a sequence of linear BVPs yielding ever smaller contributions. Those linear BVPs can then be solved straight away with RBF collocation, i.e. working in the PDE space rather than in the RBF coefficient space. %not necessarily using the same pointset, since it works at the PDE level rather than at the matrix (collocation) level. Moreover, there is no need for Jacobians, although the PDE and BCs must be analytically linearized in advance. 
This approach was used for instance in \cite{Bernal09} to solve a quasilinear PDE arising in fluid dynamics. However, we will show in Section \ref{S:Newton} that the operator-Newton method is equivalent--at least in its most straightforward version--to Newton's method for systems and thus susceptible of erratic behaviour in the event of an inadequate starting guess.

\subsection{Outline of the paper}
The remainder of the paper is organized as follows. We review the TRA in Section \ref{S:Trust}, with an eye on RBF collocation. Particular attention is paid to the so-called trust-region subproblem and three schemes for solving it are surveyed. Section \ref{S:RBF_trust} is the core of the paper, for it derives specific formulas for the Jacobian and Hessian of nonlinear RBF collocation. In Section \ref{S:Newton}, we show the equivalence between the operator-Newton approach and Newton's method for nonlinear systems of equations. Section \ref{S:PDEs} derives formulas for three important classes of nonlinear elliptic BVPs in detail, illustrating how to apply the RBF/TRA ansatz to general equations. Some comments on solvability and uniqueness of the nonlinear collocation system are made in Section \ref{S:Solvability}. Section \ref{S:Examples} reports extensive numerical experiments on four different BVPs, and Section \ref{S:Discussion} discusses them. Finally, Section \ref{S:Conclusions} concludes the paper.

\subsection{Notation}\label{SS:Notation}
\begin{itemize}
\item In the space $\mathbb{R}^d$ where the BVP $({\cal W}^{PDE},{\cal W}^{BC})$ is defined, vectors are written in bold (like ${\bf x}$ or ${\bf N}$), and operators in italics (like ${\cal W}$).
\item ${\cal L}\phi_{ij}=({\cal L}\phi)(||{\bf x}_i-{\bf x}_j||)$--like $\nabla^2\phi_{ij}$--are the entries of matrix $[{\cal L}\phi]$.
%collocated linearly transformed RBFs.
\item {\em Nodal vectors} (in ${\mathbb R}^N$) are associated to the RBF centres, such as ${\vec\alpha}$, or to a function $f({\bf x})$ evaluated over the pointset, ${\vec f}=\big(f({\bf x}_1),\ldots,f({\bf x}_N)\big)$, or to an operator ${\cal W}$ collocated on the pointset nodes, ${\vec W}$.
\item Matrices acting in the RBF coefficient space are denoted with capital letters ($J$, $H$) or like $[\phi]$, $[{\cal L}\phi]$, if they are the collocation matrix of $\phi$, ${\cal L}\phi$, etc. Finally, $diag[f]$ stands for a diagonal matrix with diagonal ${\vec f}$. 
\item With the ordering $\Xi_N=\big(\{{\bf x}_1,\ldots,{\bf x}_M\}\in\Omega\big)\cup\big(\{{\bf x}_{M+1},\ldots,{\bf x}_N\}\in\partial\Omega\big)$, $[{\cal L}\phi]_{\Omega}\in{\mathbb R}^{M\times N}$ is the upper block of the matrix $[{\cal L}\phi]$ in (\ref{F:Kansa}) and $[{\cal L}\phi]_{\partial\Omega}\in{\mathbb R}^{(N-M)\times N}$ the lower block.
\item $A>0$ ($A<0$) stands for a positive (negative) definite square matrix $A$. 
\end{itemize}

%--------------------------------------------------------------------------

\section{Overview of the trust-region algorithm}
\label{S:Trust}
%--------------------------------------------------------------------------

In this section, we discuss the TRA mainly following \cite[chapters 4 and 11]{Nocedal99}, %The TRA itself is a framework with many aspects, as well as submethods to tackle them according to the features of the main problem. Due to space limitations, we have restricted ourselves to 
focussing on those aspects which best meet the features of RBF collocation, namely: non-sparse matrices, bad conditioning, and middle-size discretizations ($N\lesssim 3000)$. Very special attention has been paid to the possibility of using the exact Hessian, which will be derived in Section \ref{S:RBF_trust}. First, a sum-of-squares scalar merit function $\mu({\vec\alpha})$ is chosen:
\begin{equation}
\label{F:Merit_function}
\mu({\vec\alpha})= {\vec W}^T{\vec W}/2=\frac{1}{2}\sum_{i=1}^NW_i^2({\vec\alpha})\geq 0.
\end{equation}  
The merit function $\mu({\vec\alpha})$ inherits the smoothness of the RBF $\phi$. Rootfinding is then recast as minimization of $\mu$:
\begin{equation}
\label{F:Minimization}
{\vec\alpha}_*= \arg \min\limits_{{\vec\alpha}\in{\mathbb R}^N} \mu({\vec\alpha}).
\end{equation}
A zero of $\mu$ is a root of the system ${\vec W}=0$, and vice versa. Moreover, a zero of $\mu$ is an absolute minimum of $\mu$. %; and, since RBFs are smooth, so are the RBF interpolant (\ref{F:RBF_Interpolant}) and $\mu({\vec\alpha})$. 
The gradient and Hessian of $\mu({\vec\alpha})$ are %the merit function (\ref{F:Merit_function}) are
\begin{equation}
\label{F:Gradient_and_Hessian}
\nabla\mu=J^T{\vec W}, \qquad H= \nabla^2 \mu= J^TJ + \sum_{i=1}^{N}W_i\nabla^2 W_i,
\end{equation}
where $J$ is the Jacobian of ${\vec W}$ (\ref{F:Jacobian_def}). Therefore, starting from an initial guess ${\vec\alpha}_0$, we seek a descending sequence $\mu({\vec\alpha}_0)>0,\mu({\vec\alpha}_1),\mu({\vec\alpha}_2),...,\mu({\vec\alpha}_{\infty})$ 
hopefully leading to the absolute minimum of $\mu$. Unfortunately, state-of-the-art minimization algorithms (not only the TRA) cannot rule out the possibility of getting trapped in a local minimum (one where $\mu>0$ and thus not a root), even if a zero of $\mu$ does exist. (The exception, as mentioned in the Introduction, are polynomial systems tackled with Groebner bases or homotopy/continuation, which pose other kind of difficulties anyway.) Moreover, minimization algorithms based on derivative information (such as the TRA) can only find {\em stationary points} (where $\nabla\mu({\vec\alpha}_{\infty})=0$), rather than minima. 
%Moreover, even if a minimum is found, it may be a local minimizer rather than a global minimum of $\mu$. Finally, if there is no root of the system in the first place, the global minimum of the merit function is strictly positive. %Summing up, depending on the problem, the RBF discretization, and the initial guess, it may happen that the sequence $\{{\vec\alpha}_k\}_{k=0}^{\infty}$ ends up in a point with singular Jacobian, a saddle point of $\mu$, or a non-zero minimum (local or strict) of the merit function, instead of a root--if there is one. 

The advantage of the TRA over other minimization algorithms is that it can deliver {\em global convergence}, which is assured convergence to {\em some} minimum of $\mu$ from {\em any} initial guess ${\vec\alpha}_0$. Precise conditions which guarantee this will be discussed later. In order to generate the iterates ${\vec\alpha}_k$, the TRA proceeds as ${\vec\alpha}_{k+1}={\vec\alpha}_k+{\vec\gamma}_k$, taking at every iteration a step ${\vec\gamma}_k$ such that
\begin{equation}
\label{F:TRS} 
{\vec\gamma}_k\approx {\vec\gamma}_{k*}:=\arg\min_{\parallel {\vec\gamma} \parallel \leq \Delta_k} \big\{ \theta_k({\vec{\gamma}})\,:=\,  \mu_k + \nabla\mu_k^T{\vec\gamma} + \frac{1}{2}{\vec\gamma}^TA_k{\vec\gamma} \big\}.
\end{equation}   
In (\ref{F:TRS}), $A_k$ is a symmetric approximation to the Hessian of $\mu_k:=\mu({\vec\alpha}_k)$, and $\Delta_k>0$ is the {\em trust-region radius}. %In practice, ${\tilde H}$ in (\ref{F:TRS}) is either ${\tilde H}=J^TJ$ (with analytical $J$ or approximated with finite differences) or ${\tilde H}=H$. 
$\theta_k$ is a second-order approximation (or {\em model}) of $\mu$ which is only {\em trusted}--and this is the hallmark of the TRA--within a distance $\Delta_k$ from the current iterate ${\vec\alpha}_k$. Problem (\ref{F:TRS})--namely, the minimization of a quadratic polynomial inside a sphere--is referred to as the {\em trust region subproblem} (TRS). Once the TRS is solved (exactly or approximately), the fidelity of the model can be assessed {\em a posteriori} by

\begin{eqnarray}
\label{F:Trust_ratio}
\rho_k= \frac{\mu({\vec\alpha}_k)-\mu({\vec\alpha}_k+{\vec\gamma}_k)}{\theta_k(0)-\theta_k({\vec\gamma}_k)}.
\end{eqnarray}

The trust-region radius can be then dynamically adjusted according to Algorithm \ref{A:TRA}. Note that a decreasing step may also be rejected if the model is deemed poor ($\rho_k\leq\eta$ for some threshold $0<\eta<1$).

\begin{algorithm}[h!]
\begin{algorithmic}
\STATE{{\bf Data:} ${\vec\alpha}_0,\Delta_{max}>0,\Delta_0\in(0,\Delta_{max})$, and $\eta\in[0,1/4)$}
\STATE{{\bf Result:} Convergence to a stationary point/minimum of $\mu$ (see Theorem \ref{Th:Convergence})}
\FOR {$k=1,2,\ldots$ until convergence}
\STATE {obtain ${\vec\gamma}_k$ by (approximately) solving (\ref{F:TRS})}
\STATE {evaluate $\rho_k$ according to (\ref{F:Trust_ratio})}
\IF {$\rho_k<1/4$}
	\STATE {$\Delta_{k+1}= \Delta_k/4$}
\ELSE
	\IF{$\rho_k>3/4$ and $\parallel {\vec\gamma}_k \parallel=\Delta_k$}
		\STATE{$\Delta_{k+1}= \min\{2\Delta_k,\Delta_{max}\}$}
	\ELSE{}
		\STATE{$\Delta_{k+1}=\Delta_k$}
	\ENDIF
\ENDIF
\IF{$\rho_k>\eta$}
	\STATE{${\vec\alpha}_{k+1}={\vec\alpha}_k+{\vec\gamma}_k$}
\ELSE{}
	\STATE{${\vec\alpha}_{k+1}={\vec\alpha}_k$}
\ENDIF
\ENDFOR
\end{algorithmic}
\caption{Trust Region Algorithm (see \cite[algorithms 4.1 and 11.5]{Nocedal99}).}
\label{A:TRA}
\end{algorithm}

There are two more important definitions. The {\em Cauchy step}, ${\vec\gamma}_C$, is the minimizer of $\theta_k$ inside the trust region along the steepest descent direction--so that $\theta_k({\vec\gamma}_C)\geq \theta_k({\vec\gamma}_*)$. It turns out to be \cite[section 4.1]{Nocedal99}:
\begin{equation}
\label{F:Cauchy_step}
{\vec\gamma}_C= -\tau_k\frac{\Delta_k}{\parallel \nabla\mu_k\parallel}\nabla\mu_k, 
\textrm{ with }\tau_k=\left\{
\begin{array}{ll}
1, & \textrm{ if $\nabla\mu_k^TA_k\nabla\mu_k\leq 0$},\\
\min\{\,\parallel \nabla_k\mu\parallel^3/(\Delta_k\nabla\mu_k^TA_k\nabla\mu_k),\,1\,\}, & \textrm{ otherwise.}
\end{array}  
\right.
\end{equation}
Note that ${\vec\gamma}_C$ involves no linear systems and yet provides some drop in $\mu$--hence it is useful to fall back on in the event of severe ill-conditioning of $A_k$. The {\em full step}, ${\vec\gamma}_F$, is the unconstrained minimizer of a quadratic polynomial: 
\begin{equation}
\label{F:Full_step}
{\vec\gamma}_F= \arg\min\limits_{{\vec\gamma}\in{\mathbb R}^N} \theta_k({\vec\gamma})= -A_k^{-1}\nabla\mu_k. 
\end{equation} 

The convergence properties of the TRA depend on the method employed to tackle the TRS, which we address next. Let us drop iteration subindex $k$ while discussing the TRS. The critical insight is that the TRS need not be solved exactly for the TRA to converge. In fact, a TRS approximation ${\vec\gamma}\approx {\vec\gamma}_*$ resulting in a drop in $\mu$ which is at least a fixed positive fraction of the drop achieved by the Cauchy step suffices for convergence \cite[theorems 4.8 and 4.9]{Nocedal99}. If $A>0$, ${\vec\gamma}_F$ is calculated. If, moreover, ${\vec\gamma}_F$ lies inside the trust region, then  ${\vec\gamma}_*={\vec\gamma}_F$. Otherwise, either $A>0$ but the full step is not feasible (i.e. $||{\vec\gamma}_F||>\Delta$), so that ${\vec\gamma}_*$ must lie on the boundary of the trust region; or $A$ is indefinite--and the full step ${\vec\gamma}_F$ may not be a minimum of $\mu$ in the first place. In both cases, the TRS--a nonlinear problem itself--must be solved iteratively. Bearing in mind the application to RBF collocation, we will consider three TRS approximations: %, depending on the computational budget and the available Hessian approximation:
\begin{enumerate}
\item {\bf Nearly exact (or "full")}. The exact solution ${\vec\gamma}_*$ of the TRS (\ref{F:TRS}) can be found following an approach due to Mor\'e and Sorensen \cite{More1983} (see also \cite[chapter 4]{Nocedal99}). The higher computational cost of this TRS method is only warranted if the full Hessian is available, so we will assume that $A=H$. A canned implementation is Matlab's \verb|trust|--albeit it is not given as an option in Matlab's \verb|fsolve|--where the full eigendecomposition of $H$ is carried out. An important case is when $\nabla\mu^T{\vec q}_1=0$ (${\vec q}_1$ is the eigenvector of $\lambda_1$, the smallest eigenvalue of $A=H$), called by Mor\'e and Sorensen the {\em hard case}, which may appear due to ill-conditioning of the RBF matrices. Alternatively, the method described in detail in \cite[section 4.2]{Nocedal99}) involves several factorizations of $H$ rather than the full eigendecomposition. Due to space limitations, the reader is referred to those papers for further details. Despite its significantly higher computational cost, this nearly exact solution will serve as a benchmark to assess the performance of the remaining two approximations.
\item The {\bf dogleg method}, involving just one matrix factorization ($A=J^TJ$).
\item {\bf 2D subspace minimization}, involving two factorizations of $A=H$.
\end{enumerate}
Conjugate-gradient-based methods for the TRS have not been included because they are mostly meant for large and sparse matrices and are especially sensitive to bad conditioning, while--to the best of our knowledge--there are not really efficient general preconditioners available for our problem.

\subsection{The dogleg method} \label{SS:Dogleg}
When the full step ${\vec\gamma}_F$ is not feasible, the minimum in the TRS is approximated by the intersection of the trust region with the "dogleg path" \cite[figure 4.3]{Nocedal99}

\begin{equation}
{\vec p}(\tau)=
\left\{
\begin{array}{ll}
\tau {\vec\gamma}_u, & \textrm{ if } 0\leq \tau \leq 1, \\
{\vec\gamma}_u + (\tau-1)({\vec\gamma}_F - {\vec\gamma}_u), & \textrm{ if } 1\leq \tau \leq 2, \\
\end{array}
\right.
\end{equation}
where
\begin{equation}
{\vec\gamma}_u= -\frac{\nabla\mu^T\nabla\mu}{\nabla\mu^TA\nabla\mu}\nabla\mu.
\end{equation}

It can be proved that $\tau\in[0,2]$ is the solution of the scalar equation
\begin{equation}
||{\vec\gamma}_u+(\tau-1)({\vec\gamma}_F - {\vec\gamma}_u)||^2=\Delta^2.
\end{equation}

See \cite[sections 4.1 and 11.2]{Nocedal99} for further details. Therefore, minimization in $N$ dimensions is replaced by minimization along the dogleg path. %The dogleg method is unexpensive (only one linear system has to be solved, namely to compute ${\vec\gamma}_u$), and especially useful when $\mu$ is convex (so that $A=H$ is guaranteed positive definite). If it is not, the TRA may well come across an indefinite Hessian, in which case a positive-definite approximation $A$ of it is needed. 
The standard choice (in fact the default in many solvers such as Matlab's {\em fsolve}, and in the remainder of this paper) is $A=J^TJ\geq 0$, which incorporates partial Hessian information without constructing $H$ in the first place, and provides second- order convergence with just first-derivative information \cite[section 11.2]{Nocedal99}. 
 
\subsection{\bf 2D subspace minimization of the TRS}\label{SS:2DSub}

Byrd, Schnabel and Schultz extended the minimization to the whole bidimensional subspace containing the dogleg path in such a way that also indefinite $A$ can be used \cite{Byrd1988,Schultz1985}. With this method, we will also understand that 
$A=H$. Practical experience shows that often, the drop in $\mu$ in the bidimensional subspace compares to that in ${\mathbb R}^N$, but at a fraction of the cost. We put together the pseudocode in Algorithm \ref{A:2dSub}.

\begin{algorithm}[htb!]
\begin{algorithmic}
\STATE{{\bf Data:} ${\nabla\mu},H,\Delta>0,tol>0$}
%\STATE{{\bf Require:} An estimate $\nu\in(-\lambda_1,-2\lambda_1]$ of the smallest eigenvalue of $H$, $\lambda_1$}
\STATE{{\bf Require:} $\nu\approx\lambda_1:=\min eig(H)$, such that $|\nu|\in(-\lambda_1,-2\lambda_1]$ if $\lambda_1<0$}
\STATE{{\bf Result:} approximate minimizer ${\vec\gamma}\approx{\vec\gamma}_*$ of the TRS in (\ref{F:TRS})}
\STATE{}
\STATE{{\bf 1}) Check definiteness of $H$}
\IF{$\nu>tol$}
	\STATE{$H>0$. Compute ${\vec\gamma}_F$ in (\ref{F:Full_step})}
	\IF{$||{\vec\gamma}_F||\leq\Delta$}
		\STATE{${\vec\gamma}={\vec{\gamma}}_F={\vec\gamma}_*$ and {\bf return}}
	\ELSE 
		\STATE{let $S_2=span[\nabla\mu,{\vec\gamma}_F]$ and {\bf goto 2})}
	\ENDIF
\ELSIF{$|\nu|<tol$}
	\STATE{$H$ is numerically singular. ${\vec\gamma}={\vec\gamma}_C$ and {\bf return}}
\ELSE %{($H$ is indefinite)}
	\STATE{$H$ is indefinite. Compute ${\vec p}=-(H-\nu I)^{-1}\nabla\mu$}
	\IF{$||{\vec p}||\leq\Delta$}
		\STATE{let $\vec q$ be a descent direction of $\mu$ and $||\vec q||=1$}
		\STATE{let ${\vec v}=-v{\vec q}$, where $v=-{\vec p}^T{\vec q}            +\sqrt{({\vec p}^T{\vec q})^2+\Delta^2-||{\vec p}||^2}$}
		\STATE{let ${\vec\gamma}={\vec p}+{\vec v}$ and {\bf return}} 
	\ELSE
	\STATE{let $S_2=span[\nabla\mu,{\vec p}]$ and {\bf goto 2})}
	\ENDIF
\ENDIF
\STATE{}
\STATE{{\bf 2}) Let $S_2=span[{\vec s}_1,{\vec s}_2]$, with $||{\vec s}_1||=||{\vec s}_2||=1$, and $P_2=[{\vec s}_1,{\vec s}_2]$ (a projector)}
\STATE{find the minimizer ${\vec\xi}_*=\sigma_1{\vec s}_1+\sigma_2{\vec s}_2=(\sigma_1,\sigma_2)^T$ of the model $\theta$ in subspace $S_2$:}
\STATE{${\vec\xi}^*=\arg\min\limits_{\sqrt{\sigma_1^2+\sigma_2^2}\leq\Delta} \mu + \nabla\mu^TP_2(\sigma_1,\sigma_2)^T + \frac{1}{2}(\sigma_1,\sigma_2)P_2^TAP_2(\sigma_1,\sigma_2)^T$}          
\STATE{${\vec\gamma}={\vec\xi}_*$ and {\bf return}}.
\end{algorithmic}
\caption{Two-dimensional subspace approximation of the TRS (2Dsub)}
\label{A:2dSub}
\end{algorithm}

%In Algorithm (\ref{A:2dSub}), note that the two-dimensional subspace, $S_2=span[{\vec s}_1,{\vec s}_2]$, depends on whether $H$ is positive definite or not. Once the vectors ${\vec s}_1$ and ${\vec s}_2$ have been computed, the two-dimensional minimization is unexpensive to carry out. Importantly, this TRS method requires a slightly larger (in absolute value) estimate $\nu$ of ${\lambda_1}$ of $H$. (Note that if $\nu=\lambda_1$ the vector ${\vec p}$ in Algorithm (\ref{A:2dSub}) could not be computed). 
Once $\nu$ is available, the descent direction ${\vec q}$ can be the eigenvector, i.e. $H{\vec q}\approx\nu{\vec q}$. In \cite{Byrd1988}, $\nu$ is meant to be efficiently estimated with the Lanczos method. % Unfortunately, our numerical experiments in Section \ref{S:Examples} show that the Lanczos method performs poorly with ill-conditioned RBF matrices.

\subsection{Convergence properties of the TRA}

The following result is a corollary from theorems 4.8 and 4.9 in \cite{Nocedal99}, which in turn were proved in \cite{Schultz1985} and \cite{More1983}, respectively. For the dogleg method, it is possible to derive conditions on $J$ rather than on $\mu$ \cite[theorems 11.8 and 11.9]{Nocedal99}.

\begin{theorem}[Convergence of the TRA with the three TRS methods in this paper]
Assume that the TRS is solved either by the nearly exact method ($A=H$), or the dogleg method ($A=J^TJ$), or by 2D subspace minimization ($A=P_2^THP_2$). Furter assume that $||A||$ is bounded above and that $\mu$ is Lipschitz continuously differentiable and bounded below on the level set $\{{\vec\alpha}\,|\,\mu({\vec\alpha})\leq \mu({\vec\alpha}_0)\}$. Then, the TRA (Algorithm \ref{A:TRA}) with constant $\eta$ converges to a stationary point% (i.e. $\nabla\mu({\vec\alpha}_{\infty})=0$)
--which may be either an absolute minimum, a local minimum, or a saddle point. 
If, additionally, that level set is compact, the TRA with nearly-exact TRS solution either converges to a minimum (local or absolute), or ${\vec\alpha}_k$ has a limit point in the level set at which second-order necessary minimality conditions hold (but not a saddle point).
\label{Th:Convergence}
\end{theorem}

The convergence of the TRA close to a non-degenerate root (i.e. where $\det J({\vec\alpha}_*)\neq 0\,$) is quadratic if $J$ is Lipschitz-continuous around $\vec\alpha_*$ and the TRS is solved exactly \cite[theorem 11.10]{Nocedal99}. Note however that, as the iterates $\vec\alpha_k$ approach the root more closely,  eventually the root will lie within the trust region, and then ${\vec\gamma}={\vec\gamma}_F$ is nearly exactly the Newton step in (\ref{F:Newton_Method}), since $J^TJ\shortrightarrow H$ because $\{W_i\approx 0\}_{i=1}^N$ in (\ref{F:Gradient_and_Hessian}). Since Newton's method converges quadratically to a non-degenerate root, so do all three TRS methods considered in this paper.

%{\bf Remark 1}. With some mild restrictions on $\mu$, and barring severe numerical instability, the TRA with nearly exact solution of the TRS and exact Hessian is guaranteed  to converge to some minimum of the merit function--but not necessarily the global minimum--for any initial guess ${\vec\alpha}_0$. This means that, even if there is a root of the underlying nonlinear system, the TRA may get trapped into a local minimum and miss it. To the author's best knowledge, this is the best that state-of-the art deterministic algorithms can do for general (i.e. non-polynomial) nonlinear systems of equations.   

\subsection{Scaling}\label{SS:Scaling}

The TRA with spherical trust regions may perform poorly when the merit function is posed with poor scaling--i.e. changes much faster along some directions than other. Linearly rescaling the coefficients vector,
%The simplest way of making up for poor scaling is to rescale the parameter space as
\begin{equation}
\label{F:Rescaling_alfa}
{\vec\alpha}':= \Gamma{\vec\alpha},\qquad (\det\Gamma\neq 0)
\end{equation}
%where $\Gamma$ is a diagonal matrix with positive elements.
%with $\det\Gamma\neq 0$. 
may make up for it. Then, 
%$\partial W/\partial\alpha_i= \partial W/\partial\alpha'_i\Gamma_{ii}$, and
\begin{equation}
\label{F:Rescaling_J}
J_{ij}=\frac{\partial W_i}{\partial\alpha_j}=\sum_{k=1}^N\frac{\partial W_i}{\partial\alpha'_k}\frac{\partial\alpha'_k}{\partial\alpha_j}=
\sum_{k=1}^NJ'_{ik}\Gamma_{kj},
\end{equation}
so that
\begin{equation}
\label{F:Rescaled_JH}
J'= J\Gamma^{-1}, \qquad A'= \Gamma^{-T}A\Gamma^{-1}. 
\end{equation} 

Notice that this may change the conditioning of $J'$ and $A'$. The simplest choice is to take $\Gamma$ diagonal with positive elements--which is equivalent to keeping the original variables and replacing the spherical trust region by an elliptical one defined by $||\Gamma\gamma_k||\leq \Delta_k$ \cite[section 4.4]{Nocedal99}. The diagonal entries $\Gamma_{ii}$ must be a reflection of the sensitivity of the merit function to changes along the $i^{th}$ coordinate axis. A reliable option is setting
\begin{equation}
\label{F:Diagonal_scaling}
%\Gamma_{ii}= \textrm{ proportional to }\frac{\partial^2\mu}{\partial\alpha_i^2}.
\Gamma_{ii}=\frac{\partial^2\mu}{\partial\alpha_i^2}.
\end{equation}

%--------------------------------------------------------------------------
\section{The trust-region algorithm for RBF collocation}
\label{S:RBF_trust}
%--------------------------------------------------------------------------

In this section, we address specific aspects of the TRA when applied to RBF collocation, including analytical formulas for $J$ and $H$. Henceforth, the combined method will be referred to as RBFTrust.%From now on, we will refer to the combined method to solve nonlinear BVPs as RBFTrust.  

\subsection{RBF Jacobian and Hessian}\label{SS:J_and_H}

The Jacobian of (\ref{F:Nonlinear_System}) is:
\begin{equation}
\label{F:Jacobian_def}
J({\vec\alpha})= 
\left[
\begin{array}{ccc}
\frac{\partial W_1}{\partial \alpha_1} & \ldots & \frac{\partial W_1}{\partial \alpha_N} \\
\vdots & \ddots & \vdots \\
\frac{\partial W_N}{\partial \alpha_1} & \ldots & \frac{\partial W_N}{\partial \alpha_N}
\end{array}
\right].
\end{equation}
Recall that $W_i$ refers to either ${\cal W}^{PDE}$ or ${\cal W}^{BC}$ depending on whether ${\bf x}_i$ is in $\Omega$ or on $\partial\Omega$, and thus $J\neq J^T$. We will assume that ${\cal W}$ is smooth and consists of $S$ functions of $u$ and its derivatives with respect to ${\bf x}$ (including the identity operator $I$). For instance, in ${\cal W}u= \nabla^2u + \sqrt{u}(\partial u/ \partial x)^2 - u$, there are three such components, namely $D_1=I,D_2=\partial u/\partial x$, and $D_3=\nabla^2$ (the order is irrelevant). Replacing $u$ by ${\tilde u}$ and applying the chain rule,
\begin{eqnarray}
\frac{\partial W_i}{\partial \alpha_j}= \sum\limits_{m=1}^S \frac{\partial W_i}{\partial D_m{\tilde u}({\bf x}_i)}\frac{\partial D_m{\tilde u}({\bf x}_i)}{\partial\alpha_j}.
\end{eqnarray}
Let us use the shorthand notation $\partial W_i/\partial D_m$ for $\partial W_i/\partial D_m{\tilde u}({\bf x}_i)$. By linearity,
\begin{eqnarray}
\frac{\partial W_i}{\partial \alpha_j}= \sum\limits_{m=1}^S \frac{\partial W_i}{\partial D_m} \frac{\partial}{\partial \alpha_j}\sum\limits_{k=1}^N \alpha_k D_m\phi_{ik}= \sum\limits_{m=1}^S \frac{\partial W_i}{\partial D_m} D_m\phi_{ij}.
\end{eqnarray}
%(Recall that $D_k\phi_{ij}= \psi({\bf x}_i)$, where $\psi({\bf x})=D_k\phi(\parallel {\bf x}-{\bf x}_j \parallel_2)$ ). 
Then, using the notation introduced in Section \ref{SS:Notation}:
\begin{equation}
\label{F:RBF_Jacobian}
J({\vec\alpha})= 
\left[
\begin{array}{ccc}
\sum\limits_{m=1}^S\frac{\partial W_1}{\partial D_m}D_m\phi_{11} & \ldots & \sum\limits_{m=1}^S\frac{\partial W_1}{\partial D_m}D_m\phi_{1N} \\
\vdots & \ddots & \vdots \\
\sum\limits_{m=1}^S\frac{\partial W_N}{\partial D_m}D_m\phi_{N1} & \ldots & \sum\limits_{m=1}^S\frac{\partial W_N}{\partial D_m}D_m\phi_{NN} \\
\end{array}
\right] = \sum\limits_{m=1}^S \,diag\Big[\frac{\partial W}{\partial D_m}\Big][D_m\phi].
\end{equation}
 
The Hessian of the merit function $\mu={\vec W}^T{\vec W}/2$ is
\begin{equation}
H= \nabla^2 \mu= J^TJ + \sum_{k=1}^{N}W_k\nabla^2 W_k, 
\end{equation}
where
\begin{eqnarray}
\nabla^2 W_k =
 \left[
\begin{array}{lcr}
\frac{\partial^2 W_k}{\partial \alpha_1^2} & \ldots & \frac{\partial^2 W_k}{\partial \alpha_1\partial \alpha_N} \\
\vdots & \ddots & \vdots \\ 
\frac{\partial^2 W_k}{\partial \alpha_N\partial \alpha_1} & \ldots & \frac{\partial^2 W_k}{\partial \alpha_N^2}\\
\end{array}
\right].
\end{eqnarray}
Notice that $H$, $J^TJ$ and $W_k\nabla^2 W_k$ are symmetric. %For future reference, let us call $J^TJ$ the {\em linear part} of $H$ (since it can be obtained from first derivatives), and $H-J^TJ$ the {\em nonlinear part}. 
Now,
\begin{eqnarray}
\frac{\partial^2 W_k}{\partial \alpha_i\alpha_j}=
\frac{\partial}{\partial \alpha_i}\Big( \sum\limits_{m=1}^S \frac{\partial W_k}{\partial D_m} D_m \phi_{kj} \Big)=
\sum\limits_{m=1}^S\sum\limits_{n=1}^S \frac{\partial^2 W_k}{\partial D_m\partial D_n}   D_n\phi_{ki} D_m\phi_{kj}.
\end{eqnarray}

%If the $L^{(0)},\ldots,L^{(M)}$ are partial derivatives, matrices $[ L^{(s)} \phi]$ must be either symmetric or antisymmetric for RBFs. For example, for the multiquadric $\phi=\sqrt{c^2+r^2}$, one has $(\partial\phi/\partial x)(\parallel {\bf x}_1 - {\bf x}_2 \parallel) = -(\partial\phi/\partial x)(\parallel {\bf x}_2 - {\bf x}_1 \parallel)/$, but $\nabla^2\phi_{12}=\nabla^2\phi_{21}$. Let us use the notation $[L^{(s)} \phi]^T= \Pi_s [L^{(s)} \phi]$, where $\Pi_s=\pm 1$. Then,

Note that $D_m\phi_{ki}=\pi_m D_m\phi_{ik}$, with $\pi_m=\pm 1$, and thus $[D_m\phi]^T=\pi_m[D_m\phi]$. For instance, $\nabla^2\phi_{ij}=\nabla^2\phi_{ji}$ but $\frac{\partial\phi}{\partial x}|_{ij}=-\frac{\partial\phi}{\partial x}|_{ji}$. Then

%\begin{eqnarray}
%\frac{\partial^2 W_k}{\partial\alpha_i\partial \alpha_j}= \sum\limits_{m=1}^S\sum\limits_{^n=1}^S \pi_nD_n\phi_{ik}\Big( \frac{\partial^2W_k}{\partial D_m\partial D_n} \Big) D_m\phi_{kj}
%\end{eqnarray}
 
\begin{equation}
\Big( \sum\limits_{k=1}^N W_k\nabla^2 W_k\Big)_{ij}=
\sum\limits_{m=1}^S\sum\limits_{n=1}^S \pi_m \sum\limits_{k=1}^N D_m\phi_{ik} \Big( W_k\frac{\partial^2 W_k}{\partial D_m\partial D_n} \Big) D_n\phi_{kj}.
\end{equation}
From (\ref{F:RBF_Jacobian}),
\begin{equation}
J^TJ= 
%\sum\limits_{m=1}^S\sum\limits_{n=1}^S [D_m\phi]^Tdiag[\frac{\partial W}{\partial D_m}]^+diag[\frac{\partial W}{\partial D_n}][D_n\phi]=
\sum\limits_{m=1}^S\sum\limits_{n=1}^S [D_m\phi]^Tdiag[\frac{\partial W}{\partial D_m}\frac{\partial W}{\partial D_n}][D_n\phi].
\end{equation}
Summing the two parts,
\begin{eqnarray}
\label{F:RBF_Hessian}
H({\vec\alpha})= \sum\limits_{m=1}^S\sum\limits_{n=1}^S [D_m\phi]^T\Big(diag[\frac{\partial W}{\partial D_m}\frac{\partial W}{\partial D_n}]+diag[W\frac{\partial^2 W}{\partial D_m\partial D_n}]\Big)[D_n\phi]= \\\nonumber
=\frac{1}{2}\sum\limits_{m=1}^S\sum\limits_{n=1}^S [D_m\phi]^Tdiag[\frac{\partial^2W^2}{\partial D_m\partial D_n}][D_n\phi].
\end{eqnarray}
Note that (\ref{F:RBF_Hessian}) is symmetric with respect to the matrices involved.\newline

{\bf Remark 2.} The matrices $[D_1\phi],\ldots,[D_S\phi]$ are filled at start and stored. At every iteration of RBFTrust, only the diagonal matrices depend on ${\vec\alpha}_k$ and have to be recalculated. Thus the dogleg method--where only $J$ is explicitly computed--involves only matrix-vector multiplications, while setting $A=H$ takes $S(S+1)/2$ extra matrix multiplications. 

\subsection{Elimination of linear equations}

Often, the system ${\vec W}(\vec\alpha)=0$ will contain linear equations, such as those representing the collocation of Dirichlet and other linear BCs. Another source of linear equations in the system is the enrichment of the RBF interpolant (see for instance \cite{Bernal_Bilap,Bernal_Newt}, and Example III in Section \ref{S:Examples}) with $n$ special functions $h_k$: 

%(\ref{}) with augmenting functions, for instance tailored to capture boundary singularities (see Example IV in Section \ref{S:Examples}):

\begin{equation}
\label{F:Enriched_RBF_Interpolant}
{\tilde u}= \sum_{j=1}^N \alpha_j \phi(\parallel {\bf x}-{\bf x}_j\parallel_2) + \sum_{k=1}^n h_k({\bf x}).
\end{equation}
In this case, the system must be augmented with $n$ ancillary equations in order to keep it square: 
\begin{equation}
\label{F:Complementary_equations}
\alpha_1h_k({\bf x}_1)+\ldots+\alpha_Nh_k({\bf x}_N)=0,\qquad k=1,\ldots,n.
\end{equation}

Whenever there are $m\leq N$ linearly independent equations, $m$ degrees of freedom can be eliminated, and the minimum for $\mu$ sought in a shrunken $(N-m)$-dimensional space. An elimination method which is optimally stable %and well-suited to RBF collocation 
is described in \cite[section 15.2]{Nocedal99}.
Assume that the linear block in (\ref{F:Nonlinear_System}) is given by $B{\vec\alpha}={\vec b}$, with $B\in{\mathbb R}^{m\times N}$. Consider the QR decomposition 
\begin{equation}
\label{F:QR}
B^T\Pi= [Q_1 Q_2]
\left[\begin{array}{c}
R \\ 0
\end{array}\right],
\end{equation} 
where $\Pi$ is an $m\times m$ permutation matrix, $Q_B=[Q_1 Q_2]$ is orthonormal, $Q_1\in{\mathbb R}^{N\times m}$ and $Q_2\in{\mathbb R}^{N\times (N-m)}$ are made up of orthonormal columns and $R\in{\mathbb R}^{m\times m}$ is upper triangular and nonsingular because $rank(B)=m$. Let ${\vec\alpha}= Q_1{\vec\nu} + Q_2{\vec\beta}$ and insert it into $B{\vec\alpha}={\vec b}$, yielding the optimally stable decomposition of ${\vec\alpha}$

\begin{equation}
\label{F:alfa_beta}
{\vec\alpha}= Q_1R^{-T}\Pi^T{\vec b} + Q_2{\vec\beta}.
\end{equation}

In any elimination of degrees of freedom like $\vec{\alpha}=\vec{v} + Z\vec{\beta}$, where $\vec{v}$ is a constant vector and $\vec{\beta}\in{\mathbb R}^{N-m}$, the columns of $Z$ represent a basis of the null space of $B$, since $BZ=0$. (It is easy to prove that $BQ_2=0$). Whether $Z=Q_2$ or a different basis of $null(B)$, Jacobian entries are transformed as in (\ref{F:Rescaling_J}) and
%\begin{equation}
%\frac{dW_i}{d\beta_j}=\sum_{k=1}^N\frac{dW_i}{d\alpha_k}
%\frac{\partial\alpha_k}{\partial\beta_j}= \sum_{k=1}^N\frac{dW_i}{d\alpha_k}Z_{kj}=[J(\vec{\alpha})Z]_{ij},
%\end{equation}
%and thus
\begin{equation}
\label{F:small_J_and_H}
J(\vec{\beta})= J(\vec{\alpha})Z, \qquad H(\vec{\beta})=Z^TH(\vec{\alpha})Z.
\end{equation}

In (\ref{F:small_J_and_H}), $J(\vec{\alpha})$ and $H(\vec{\alpha})$ are (\ref{F:RBF_Jacobian}) and (\ref{F:RBF_Hessian}), respectively. (Or $J'$ and $H'$ like in (\ref{F:Rescaled_JH}) if rescaling like (\ref{F:Rescaling_alfa}) has been included.)

 %With the specific, optimally stable choice just described, note that indeed $BZ=...=0$.    

%--------------------------------------------------------------------------
\section{The operator-Newton approach}
\label{S:Newton}
%--------------------------------------------------------------------------

Let us define the linearization of the nonlinear operator ${\cal W}$ %acting between Banach spaces
around a function $u({\bf x})$ as the linear operator ${\cal L}_u$ such that
\begin{equation}
\label{F:Frechet1}
\lim\limits_{\parallel v\parallel\rightarrow 0}\frac{{\cal W}(u+v)-{\cal W}u-{\cal L}_uv}{\parallel v \parallel} =0.
\end{equation} 
If that limit exists, ${\cal L}_u$ is unique and is called the Fr\'echet derivative of ${\cal W}$ around $u$ \cite{Dierkes2002}. Definition (\ref{F:Frechet1}) is non-constructive. For our purposes, let us assume that ${\cal L}_u$ exists and that all operators are regular enough so that
\begin{equation}
\label{F:Frechet2}
{\cal W}(u+v)= {\cal W}(u) + {\cal L}_uv + {\cal O}(\parallel v\parallel^2). %\textrm{  if $\parallel v\parallel<<1$,}
\end{equation}     
%Let us change the notation so that ${\cal W}u=0$ stands for the full BVP in the sense of (\ref{F:Nonlinear_Elliptic_BVP}). Then, the linear operators 
Let ${\cal L}_u^{PDE}$ and ${\cal L}_u^{BC}$ be the linearization of ${\cal W}^{PDE}$ and ${\cal W}^{BC}$, respectively. The operator-Newton method for nonlinear elliptic BVPs is Algorithm \ref{A:opNewton} below. 

\begin{algorithm}[h!]
\begin{algorithmic}
\STATE{{\bf Data:} initial guess $u_0({\bf x})$, the linearized operators ${\cal L}_u=\{{\cal L}_u^{PDE},{\cal L}_u^{BC}\}$}
%\STATE{{\bf Require:} the linearized operator ${\cal L}_u=\{{\cal L}_u^{PDE},{\cal L}_u^{BC}\}$ }
\FOR{$k=1,2,...$ until $R_k=\{R^{PDE}_k,R^{BC}_k\}$ stagnates}
\STATE{\begin{align*}
&\textrm{{\bf A3.1.} Compute the residuals} \left\{
\begin{array}{ll} 
R_k^{PDE}=-{\cal W}^{PDE}u_{k-1} & \textrm{ if }{\bf x}\in\Omega,\\ & \\
R_k^{BC}= -{\cal W}^{BC}u_{k-1} &  \textrm{ if }{\bf x}\in\partial\Omega\\
\end{array}\right.\\
&\textrm{{\bf A3.2.} Solve the BVP} \left\{
\begin{array}{ll} 
{\cal L}^{PDE}_{u_k}v_k=R^{PDE}_k & \textrm{ if }{\bf x}\in\Omega,\\ & \\
{\cal L}^{BC}_{u_k}v_k=R^{BC}_k &   \textrm{ if }{\bf x}\in\partial\Omega\\
\end{array}\right.\\
&\textrm{{\bf A3.3.} Update $u_k=u_{k-1}+v_k $}
\end{align*}}
\ENDFOR
\end{algorithmic}
\caption{Operator-Newton method without smoothing.}
\label{A:opNewton}
\end{algorithm}

\subsection{Equivalence to Newton's method} 
Note that Algorithm \ref{A:opNewton} is independent of the discretization. Let us assume that every iteration of it is implemented on a pointset $\Xi_N$ using the RBF $\phi$, i.e. ${\tilde u}_k=[\phi]{\vec \alpha}_k$ and ${\tilde v}_k=[\phi]{\vec\gamma}_k$.
%$u\approx {\tilde u}=\sum_{j=1}^N\alpha_j\phi(\parallel{\bf x}_j-{\bf x}\parallel)$ and $v\approx {\tilde v}_k=\sum_{j=1}^N\gamma_j\phi(\parallel{\bf x}_j-{\bf x}\parallel)$. 
%and let ${\tilde u}_k$ and ${\tilde v}_k$ be the RBF approximations to $u_k$ and $v_k$ with nodal coefficients ${\vec \alpha}_k$ and ${\vec\gamma}_k$, respectively. 
Then, the updating step ({\bf A3.3}) is equivalent to ${\vec\alpha}_{k+1}={\vec\alpha}_k+{\vec\gamma}_{k+1}$. 
%By the chain rule, the collocated version of ({\bf A.3.2}) yields: 
%\begin{equation}
%{\cal L}_{u_k} v({\bf x}_i) %\approx L_W {\tilde v}({\bf x}_i)
%\approx \Big( {\cal L}_{u_k}\sum\limits_{j=1}^N \gamma_j \phi(\parallel{\bf x}-{\bf x}_j\parallel\Big)({\bf x}_i)= \sum\limits_{j=1}^N \gamma_jL_k\phi_{ij}, 
%\end{equation}
%where $L_k\phi_{ij}={\cal L}_{u_k}\phi(||{\bf x}-{\bf x}_j||)({\bf x}_i)$. 
The matrix version of an iteration of Algorithm \ref{A:opNewton} reads $[{\cal L}_k\phi]{\vec\gamma}_k=-{\vec W}_k$, 
--i.e. exactly as in (\ref{F:Kansa}) but replacing ${\vec\alpha}$ by ${\vec\gamma}_k$, $({\cal L}^{PDE},{\cal L}^{BC})$ by $({\cal L}_k^{PDE},{\cal L}_k^{BC})$, and the right-hand side by $-{\vec W}_k$.

Inverting that system for ${\vec\gamma}_k$, the matrix form of the updating step ({\bf A3.3}) in Algorithm \ref{A:opNewton} gives ${\vec\alpha}_{k+1}={\vec\alpha}_k-[{\cal L}_k\phi]^{-1}{\vec\gamma}_{k+1}$, which is Newton's method if $[{\cal L}_k\phi]=J_k$. %It is well known that in finite-dimensional spaces (like the space spanned by the RBFs in $\Xi_N$, here) the Fr\'echet derivative is the Jacobian. 
To justify that this is indeed the case, consider the limit
\begin{equation}
\lim\limits_{\parallel v \parallel\to 0} \frac{{\cal W}(u+v)-{\cal W}u}{\parallel v\parallel}=
\lim\limits_{\parallel v \parallel\to 0} \frac{{\cal L}_uv}{\parallel v\parallel},
\end{equation}
according to (\ref{F:Frechet1}) and (\ref{F:Frechet2}). Let ${\vec\gamma}_k=(\gamma_k^{(1)},\ldots,\gamma_k^{(N)})^T$. Substituting ${\tilde u}({\vec\alpha}_k)$ and ${\tilde v}({\vec\gamma}_k)$ for $u$ and $v$ and evaluating on ${\bf x}_i\in\Xi_N$,
\begin{equation}
\label{F:aux1}
\lim\limits_{\parallel {\vec\gamma}_k \parallel\to 0} \frac{{\cal W}_i({\vec\alpha}_k+{\vec\gamma}_k)-{\cal W}_i({\vec\alpha}_k)}{\parallel {\vec\gamma}_k\parallel}=
\lim\limits_{\parallel {\vec\gamma}_k \parallel\to 0} \frac{\sum\limits_{n=1}^N\gamma_k^{(n)} {\cal L}_k\phi_{in}}{\parallel {\vec\gamma}\parallel}.
\end{equation}
%Now, let us particularize to $\gamma=(0,\ldots,1,0,\ldots,0)$, where the $1$ is on the $j^{th}$ position. Then, the left-hand side of (\ref{F:aux1}) is $\partial W_i/\partial\alpha_j=J_{ij}$, and
Now, let us particularize to the direction $\gamma_k^{(n)}=\delta_{jn}$, where $\delta_{jn}$ is Kronecker's delta. Then, the left-hand side of (\ref{F:aux1}) is $\partial W_i/\partial\alpha_j=J_{ij}$, and
\begin{equation}
J_{ij}= L_k\phi_{ij}. %\box
\end{equation} 
Thus, the collocated version of Algorithm \ref{A:opNewton} is not globally convergent.\newline % , and in general cannot be expected to converge to the solution of the nonlinear equation from any initial hint of the solution.\newline

{\bf Remark 3.} The original operator-Newton algorithm by Fasshauer includes %the possibility of using nested pointsets as well as 
an additional step for residual smoothing, and therefore our result here does not pertain to that case. In fact, his research shows that smoothing may be critical for performance, although difficult to implement in practice \cite{Fasshauer00,Bernal09}.

%--------------------------------------------------------------------------
\section{Explicit formulas for prototypical PDEs}
\label{S:PDEs}
%--------------------------------------------------------------------------
In this section, we derive explicit formulas for the Jacobian and Hessian of three relevant classes of nonlinear elliptic differential operators, as well as for the linearized operator required by the operator-Newton method. Those formulas will be used later in Examples I-IV in Section \ref{S:Examples}.

\subsection{Semilinear equations}\label{SS:Semilinear}

In this case, the nonlinearity involves $u$ but none of its derivatives. We consider only the following  equation, taken from \cite{Xu93}:
\begin{eqnarray}
\label{F:PDE_cubic}
\nabla^2 u - u^3= f \textrm{ if }{\bf x}\in\Omega,\qquad
u= g \textrm{ if }{\bf x}\in\partial\Omega.
\end{eqnarray}

The Fr\'echet derivatives for the operator-Newton method are 
\begin{equation}
{\cal L}^{PDE}=\nabla^2-3u^2I,\qquad {\cal L}^{BC}=I.
\end{equation} 

For RBFTrust, since the BCs are linear, $Z$ is obtained (along with $\Pi$ and $Q_1$) from the QR decomposition of the block with the BCs, arranged as in (\ref{F:QR}):
\begin{equation}
\label{F:Z_Dirichlet}
[\phi]_{\partial\Omega}^T\Pi= [Q_1 Z]\left[\begin{array}{c}R\\0\end{array}\right].
\end{equation}
Recall from (\ref{F:alfa_beta}) that, after finding the solution ${\vec\beta}_{\infty}$ in the shrunken space ${\mathbb R}^M$, the coefficients of the RBF interpolant are transformed according to
\begin{eqnarray}
{\vec\alpha}({\vec\beta})= Q_1R^{-T}\Pi^T[g({\bf x}_{M+1}),\ldots,g({\bf x}_N)]^T + Z{\vec\beta}.
\end{eqnarray}
The Jacobian and Hessian in the shrunken space are given by
\begin{equation}
\label{F:J_semilinear}
J({\vec\beta})=\Bigg(\, [\nabla^2 \phi]_{\Omega} -3diag[u^2({\vec\beta})]_{\Omega}[\phi]_{\Omega} \Bigg)\,Z,
\end{equation}
\begin{eqnarray}
\label{F:H_semilinear}
H({\vec\beta})=Z^T\,\Bigg(\, [\nabla^2\phi]^2_{\Omega} + [\phi]_{\Omega}diag[15u^4({\vec\beta})-6u({\vec\beta})\nabla^2u({\vec\beta})]_{\Omega}[\phi]_{\Omega}-\nonumber\\ 
-\big([\phi]_{\Omega}+[\nabla^2\phi]_{\Omega}\big)diag[3u^2({\vec\beta})]_{\Omega}\big([\phi]_{\Omega}+[\nabla^2\phi]_{\Omega}\big)\, \Bigg)\,Z.
\end{eqnarray}

The nodal values for the diagonal matrices are picked from ${\vec u}({\vec\beta})=[\phi]{\vec\alpha}({\vec\beta})$.

\subsection{Quasilinear equations}\label{SS:Quasilinear}

Let us consider the quasilinear operator
\begin{equation}
\label{F:Quasilinear_PDE}
{\cal W}^{PDE}u= \nabla\cdot\big( G(u)\nabla u \big) - f,
\end{equation}
where %$\nabla\cdot$ is the divergence operator and 
$G(t):[0,\infty)\rightarrow{\mathbb R}$ is smooth and $t$ actually stands for $|\nabla u|$,
i.e. $G(u) = G(|\nabla u|)$. %(However, for ease of notation, we will stick to the notation G(u).) 
This class includes:
\begin{equation}
\left\{
\begin{array}{l}
\textrm{ the operator $G(t)=1/\sqrt{1+t^2}$ (see Example $II$), and}\\
\textrm{ the {\em p-Laplacian}, $G(t)=t^{p-2},\,(p\geq 2)$ (see Example $III$). }
\end{array}
\right.\nonumber
\end{equation}

First, we will produce a linearization of the PDE suitable for the operator- Newton using a "small increment" argument. Consider the gradient modulus: 
\begin{equation}
|\nabla(u+v)|= 
%\sqrt{(\nabla u + \nabla v \big)^2}= 
\sqrt{ |\nabla u|^2 + |\nabla v|^2 + 2(\nabla u \cdot \nabla v) }=
|\nabla u|\sqrt{1 + 2\frac{\nabla u \cdot \nabla v}{|\nabla u|^2} + \big(\frac{|\nabla v|}{|\nabla u|}\big)^2}.
\end{equation}

Under the condition $|\nabla v|/|\nabla u|<<1$, a Taylor approximation yields

\begin{equation}
|\nabla(u+v)|= |\nabla u|\,\Big( 1 + \frac{\nabla u \cdot \nabla v}{|\nabla u|^2} + 
\frac{1}{2}\big(\frac{|\nabla v|}{|\nabla u|}\big)^2 \Big) \approx
|\nabla u| + \frac{\nabla u \cdot \nabla v}{|\nabla u|} + {\cal O}(|\nabla v|/|\nabla u|)^2.
\end{equation}

Using the definition (\ref{F:Frechet2}), it is clear that ${\cal L}^{(|\nabla|)}u= \frac{\nabla u \cdot \nabla}{|\nabla u|}$. By the chain rule,

\begin{equation}
{\cal L}^{(G)}u= G'(u)\frac{\nabla u \cdot \nabla}{|\nabla u|}
\end{equation}  

and therefore, to first order in $|\nabla v|/|\nabla u|$,

\begin{equation}
G(u+v)\approx Gu + G'(u)\frac{\nabla u \cdot \nabla v}{|\nabla u|}.
\end{equation}

%Let us now linearize the divergence in (\ref{F:Quasilinear_PDE}):
with $G,G'$ and $G''$ evaluated at $|\nabla u|$. Coming back to (\ref{F:Quasilinear_PDE}),
\begin{eqnarray}
\label{F:Divergence_of_G}
\nabla\cdot[G(u+v)\nabla(u+v)]= 
\nabla\cdot[G\nabla u] + \nabla\cdot[G\nabla v] + \nonumber\\
+ \nabla\cdot[G' \frac{\nabla u \cdot \nabla v }{|\nabla u|} \nabla u] 
+ \nabla\cdot[G' \frac{\nabla u \cdot \nabla v}{|\nabla u|} \nabla v] + O(|\nabla v|/|\nabla u|^2).
\end{eqnarray}

The first three terms on the rhs can be simplified to: 
\begin{eqnarray}
\label{F:3Terms}
\nabla\cdot[G\nabla u]= G\nabla^2 u +  G' \frac{\Delta_{\infty}u}{|\nabla u|}\\
\nabla\cdot[G\nabla v]= G\nabla^2 v + G'\frac{\nabla v \cdot (\nabla u\cdot \nabla)\nabla u}{|\nabla u|} \\
\nabla\cdot[G' \frac{\nabla u \cdot \nabla v }{|\nabla u|} \nabla u]=
G'\frac{\nabla u \cdot (\nabla u\cdot\nabla)\nabla v}{|\nabla u|} + \nonumber\\
+\Big( G''\frac{\Delta_{\infty}u}{|\nabla u|^2} - G'\frac{\Delta_{\infty}u}{|\nabla u|^3} + G'\frac{\nabla^2 u}{|\nabla u|} \Big)\nabla u \cdot \nabla v +
G'\frac{\nabla v \cdot (\nabla u\cdot \nabla)\nabla u}{|\nabla u|},
\end{eqnarray}

where we %have used the identities .... and 
have introduced the infinity Laplacian in ${\mathbb R}^d$, $\Delta_\infty=\sum_{i,j=1}^d\frac{\partial ^2}{\partial x_i\partial x_j}\frac{\partial}{\partial x_i}\frac{\partial}{\partial x_j}$. 
Regarding the rightmost term of (\ref{F:Divergence_of_G}), which is nonlinear in $v$, let
\begin{equation}
{\bf V}:= G'\frac{\nabla u \cdot \nabla v}{|\nabla u|}\nabla v.
\end{equation}
In order that $|{\bf V}| \leq |G'||\nabla v|^2 \leq \Big( \frac{|\nabla v|}{|\nabla u|} \Big)^2 \approx 0$, it is sufficient that 
\begin{equation}
\label{F:Linearization_condition}
|G'|\leq \frac{1}{|\nabla u|^2},\qquad\textrm{ or }\qquad| G'(t)|\leq 1/t^2.
\end{equation}
Thus, if (\ref{F:Linearization_condition}) holds, the surviving nonlinear term in (\ref{F:Divergence_of_G}) can be dropped. In the operator-Newton method, the correction $v$ around $u$ obeys the linear PDE:
\begin{eqnarray}
\label{F:PDE_for_v}
G\nabla^2 v + G'\frac{\nabla u}{|\nabla u|}\cdot(\nabla u \cdot \nabla)\nabla v + 
\frac{ G''|\nabla u|\Delta_{\infty}u + G'|\nabla u|^2\nabla^2 u - G'\Delta_{\infty}u}{|\nabla u|^3}\nabla u \cdot \nabla v + \nonumber\\
+ 2\frac{G'}{|\nabla u|}\nabla v \cdot (\nabla u \cdot \nabla)\nabla u = R^{PDE}
\textrm{ if ${\bf x}\in\Omega$, }
%\\{\cal L}^{B}v= R^{BC}
\end{eqnarray}

where $R^{PDE}= f-div[G(u)\nabla u]$ 
%\begin{equation}
%R^{PDE}= f-div[G(u)\nabla u] %, \qquad\qquad R^{BC}= -B_uv
%\end{equation}
%are the nonlinear residual of the current iterate $u$ to the PDE and the BCs, respectively.
is the residual to (\ref{F:Quasilinear_PDE}). In $d=2$, (\ref{F:PDE_for_v}) is

%Let $u_0({\bf x})$ be the initial guess and $v_k(\bf x)$ the $k^{th}$ correction in Algorithm \ref{A:opNewton}. Recall that they are represented by the operator ${\cal W}^{BC}$ in (\ref{F:Nonlinear_Elliptic_BVP}).
%\begin{equation}
%\label{F:BC_linearization}
%{\cal L}^{(B)}v= -Bu \textrm{  if ${\bf x}\in\partial\Omega$,}
%\end{equation} 
%where $B$ is the nonlinear boundary operator which complements the PDE (\ref{F:Two_quasilinear_PDEs}). 

%For better numerical stability it may be preferable to solve
%\begin{eqnarray}
%G|\nabla u|^3 \nabla^2 v + G'|\nabla u|^2\nabla u\cdot(\nabla u \cdot \nabla)\nabla v + \nonumber\\
%\Big(G''|\nabla u|\Delta_{\infty}u + G'|\nabla u|^2\nabla^2 u - G'\Delta_{\infty}u\Big) 
%\nabla u \cdot \nabla v + \nonumber\\
%2G'|\nabla u|^2 \nabla v \cdot (\nabla u \cdot \nabla)\nabla u = |\nabla u|^3 R.
%\end{eqnarray}

\begin{equation}
A\frac{\partial^2 v}{\partial x^2} +
B\frac{\partial^2 v}{\partial x\partial y} +
C\frac{\partial^2 v}{\partial y^2} +
D\frac{\partial v}{\partial x} +
E\frac{\partial v}{\partial y} = |\nabla u|^3R^{PDE}
\end{equation}

where:

\begin{equation}
A(x,y)= G|\nabla u|^3 + G'|\nabla u|^2u_x^2, 
\end{equation}

\begin{equation}
B(x,y)= 2G'|\nabla u|^2u_xu_y,
\end{equation}

\begin{equation}
C(x,y)= G|\nabla u|^3 + G'|\nabla u|^2u_y^2, 
\end{equation}

\begin{eqnarray}
D(x,y)= \Big( G''|\nabla u|\Delta_{\infty}u + G'|\nabla u|^2(3u_{xx}+u_{yy}) - G'\Delta_{\infty}u \Big)u_x + 2G'|\nabla u|^2u_{xy}u_y, 
%= \nonumber\\ \Big( G''|\nabla u|^2\Delta_{\infty}u + G'(u_x^2\nabla^2 u + 2|\nabla u|^2u_{xx})  \Big)u_x + 2G'u_{xy}u_y^3 MAL CALCULADO!
\end{eqnarray}

\begin{eqnarray}
E(x,y)= \Big( G''|\nabla u|\Delta_{\infty}u + G'|\nabla u|^2(3u_{yy}+u_{xx}) - G'\Delta_{\infty}u \Big)u_y + 2G'|\nabla u|^2u_{xy}u_x. 
%\nonumber\\ \Big( G''|\nabla u|^2\Delta_{\infty}u + G'(u_y^2\nabla^2 u + 2|\nabla u|^2u_{yy})  \Big)u_y + 2G'u_{xy}u_x^3 MAL CALCULADO!
\end{eqnarray}

The linearized PDE is elliptic as long as $B^2-4AC<0$, i.e. $1+t\frac{G'}{G}>0$ or
\begin{equation}
\label{F:Ellipticity_condition}
G'/G>-1/t %\frac{G'}{G}>-\frac{1}{t}
\end{equation}

Let us check the linearization (\ref{F:Linearization_condition}) and ellipticity (\ref{F:Ellipticity_condition}) conditions for the PDEs (\ref{F:Quasilinear_PDE}). For the least area operator, $G'(t)=-tG^3(t)$ and $0\leq G\leq 1$, so that both (\ref{F:Ellipticity_condition}) and (\ref{F:Linearization_condition}) hold for any $|\nabla u|$. Regarding the p-Laplace operator, (\ref{F:Ellipticity_condition}) holds, but (\ref{F:Linearization_condition}) leads to $(p-2)t^{p-1}\leq 1$%, which depends on $p$ (recall that $p\geq 2$)
. For instance, for $p=2.6$ as in Example III, the linearization breaks down if $|\nabla u|>(p-2)^{1/(1-p)}=1.37$.\newline

{\bf Remark 4.} The derivation of the linearized operator in terms of a "small increment" argument, like above, has the advantage that it sheds light on why the operator-Newton method breaks down in the presence of high gradients of the solution: it neglects contributions which are too large to be dropped, namely the nonlinear term in $v$ on the right-hand side of (\ref{F:Divergence_of_G}).\newline  

We write down now the formulas for $J$ and $H$ needed in RBFTrust. Expanding the divergence before linearization, (\ref{F:Quasilinear_PDE}) reads:
\begin{equation}
\label{F:Quasilinear_PDE2}
G\nabla^2 u + \frac{G'\Delta_{\infty}u}{|\nabla u|}= f
\end{equation}

For Plateau's problem (Example II) $f=0$, $d=2$ and, since $1+|\nabla u|^2\neq 0$, 
\begin{equation}
W^{PDE}= (1+|\nabla u|^2)\nabla^2 u - \Delta_{\infty}u,\qquad W^{BC}=u
\end{equation}
so that the non-zero entries of the RBF Jacobian (\ref{F:RBF_Jacobian}) can be filled up with: 
\begin{eqnarray}
\label{F:Jacobian_of_Plateau}
\frac{\partial W^{PDE}}{\partial u_{xx}}= 1+u_y^2,\,\,\,
\frac{\partial W^{PDE}}{\partial u_{xy}}= -2u_x u_y,\,\,\,
\frac{\partial W^{PDE}}{\partial u_{yy}}= 1+u_x^2, \nonumber\\
\frac{\partial W^{PDE}}{\partial u_x}= 2u_x u_{yy} - 2u_y u_{xy},\qquad
\frac{\partial W^{PDE}}{\partial u_y}= 2u_y u_{xx} - 2u_x u_{xy},\qquad
\frac{\partial W^{BC}}{\partial u}= 1.%\nonumber\\
%\frac{\partial W^{BC}}{\partial u_{xx}}= 
%\frac{\partial W^{BC}}{\partial u_{xy}}= 
%\frac{\partial W^{BC}}{\partial u_{yy}}= 
%\frac{\partial W^{BC}}{\partial u_x}=
%\frac{\partial W^{BC}}{\partial u_y}= 
%\frac{\partial W^{PDE}}{\partial u}= 0.  \nonumber\\
\end{eqnarray}
The non-zero second derivatives needed for the RBF Hessian are:
\begin{eqnarray}
\label{F:Hessian_of_Plateau}
\frac{\partial^2 W^{PDE}}{\partial u_x\partial u_{xy}}= -2u_y, \qquad
\frac{\partial^2 W^{PDE}}{\partial u_y\partial u_{xy}}= -2u_x, \qquad
\frac{\partial^2 W^{PDE}}{\partial u_x\partial u_{yy}}= 2u_x, \qquad
\frac{\partial^2 W^{PDE}}{\partial u_y\partial u_{xx}}= 2u_y, \nonumber\\
\frac{\partial^2 W^{PDE}}{\partial u_x\partial u_y}= -2u_{xy},\qquad
\frac{\partial^2 W^{PDE}}{\partial u_x^2}= 2u_{yy},\qquad
\frac{\partial^2 W^{PDE}}{\partial u_y^2}= 2u_{xx}.
\end{eqnarray}

For the Hele-Shaw equation, $f=0$ also, but the BCs either of Neumann or Dirichlet kind, depending on the point on the boundary (see Example III): 
\begin{equation}
{\cal W}^{PDE}= |\nabla u|^{p-2}\Big(\,\nabla^2 u + (p-2)|\nabla u|^{p-2}\Delta_{\infty}u\,\Big),\qquad 
{\cal W}^{BC}=\left\{
\begin{array}{l} {\cal W}^{N}={\bf N}\cdot\nabla u \\ {\cal W}^{D}=u-g \end{array}\right.
\end{equation}

In $d=2$, ${\bf N}=(N_x,N_y)$. The first and second derivatives of ${\cal W}^{PDE}$ with respect to $u_{xx},u_{xy},u_{yy},u_x$ and $u_y$ are found as before. (Since they are rather lengthy we do not write them down.) Since the BCs are linear, only the first derivatives of ${\cal W}^{BC}$ are nonzero, namely

\begin{equation}
\frac{\partial W^{D}}{\partial u}= 1,\qquad
\frac{\partial W^{N}}{\partial u_x}= N_x,\qquad 
\frac{\partial W^{N}}{\partial u_y}= N_y.
\end{equation}

When constructing the Jacobian and Hessian, notice that $[\phi_x]$ and $[\phi_y]$ are antisymmetric. The linear block entering the QR decomposition (\ref{F:QR}) includes not only the BCs (both Neumann and Dirichlet), but also the complementary equations (\ref{F:Complementary_equations}) for the Motz functions (\ref{F:Motz}) if they are incorporated into the RBF interpolant, as in Example III.

\subsection{Fully nonlinear operator}\label{SS:Fully_nolinear}

Here, the nonlinearity involves the highest-order derivatives. In this class, we consider the Monge-Amp\`ere operator in ${\mathbb R}^d$,
\begin{equation}
\label{F:Monge-Ampere}
{\cal W}^{PDE}u= \det{D_{d}^2u},
\end{equation}
where $D_{d}^2u$ is the Hessian matrix of $u:{\mathbb R}^d\shortrightarrow{\mathbb R}$. In $d=2$ and $d=3$,
\begin{equation}
D_2^2u= u_{xx}u_{yy}-u_{xy}^2, \qquad
D_3^2u= \det 
\left[
\begin{array}{lll}
\frac{\partial^2 u}{\partial x^2} & \frac{\partial^2 u}{\partial x\partial y} & \frac{\partial^2 u}{\partial x\partial z} \\
\frac{\partial^2 u}{\partial y\partial x} &\frac{\partial^2 u}{\partial y^2} &  \frac{\partial^2 u}{\partial y\partial z} \\ 
\frac{\partial^2 u}{\partial z\partial x} & \frac{\partial^2 u}{\partial z\partial y} &  \frac{\partial^2 u}{\partial z^2}
\end{array}
\right]
\end{equation}

This PDE gives rise to a polynomial system of collocation equations and is simpler to linearize. We assume Dirichlet BCs and write down the formulas for $d=2$. For the operator-Newton method,
\begin{equation}
{\cal L}^{PDE}= u_{yy}\frac{\partial^2}{\partial x^2} + u_{xx}\frac{\partial^2}{\partial y^2} - 2u_{xy}\frac{\partial^2}{\partial x\partial y}.
\end{equation}

For RBFTrust, $Z,Q_1$ and $\Pi$ are given by (\ref{F:Z_Dirichlet}), and the Jacobian and Hessian are:
\begin{equation}
J({\vec\beta})= \Big(\, diag[u_{yy}({\vec\beta})]_{\Omega}[\phi_{xx}]_{\Omega} + diag[u_{xx}({\vec\beta})]_{\Omega}[\phi_{xx}]_{\Omega} - 2diag[u_{xy}({\vec\beta})]_{\Omega}[\phi_{xy}]_{\Omega} \,\Big)Z
\end{equation}
\begin{eqnarray}
H({\vec\beta})= \frac{1}{2}Z^T\Bigg(\, 
  [\phi_{xx}]_{\Omega} diag[2u_{yy}^2({\vec\beta})]_{\Omega} [\phi_{xx}]_{\Omega} 
+ [\phi_{yy}]_{\Omega} diag[2u_{xx}^2({\vec\beta})]_{\Omega} [\phi_{yy}]_{\Omega}		+ \\ \nonumber																		  
  [\phi_{xy}]_{\Omega} diag[12u_{xy}^2({\vec\beta})-4u_{xx}({\vec\beta})u_{yy}({\vec\beta})]_{\Omega} [\phi_{xy}]_{\Omega} +					\\ \nonumber
\Big([\phi_{xx}]_{\Omega}+[\phi_{yy}]_{\Omega}\Big)diag[4u_{xx}({\vec\beta})u_{yy}({\vec\beta})-2u_{xy}^2({\vec\beta})]_{\Omega}\Big([\phi_{xx}]_{\Omega}+[\phi_{yy}]_{\Omega} \Big)+\\ \nonumber
\Big([\phi_{xx}]_{\Omega}+[\phi_{xy}]_{\Omega}\Big)diag[-4u_{yy}({\vec\beta})u_{xy}({\vec\beta}) ]_{\Omega}\Big([\phi_{xx}]_{\Omega}+[\phi_{xy}]_{\Omega} \Big) +\\ \nonumber
\Big([\phi_{yy}]_{\Omega}+[\phi_{xy}]_{\Omega}\Big)diag[-4u_{xx}({\vec\beta})u_{xy({\vec\beta})} ]_{\Omega}\Big([\phi_{yy}]_{\Omega}+[\phi_{xy}]_{\Omega} \Big) \Bigg)\,Z.
\end{eqnarray}

%--------------------------------------------------------------------------
\section{Remarks on solvability and uniqueness}
\label{S:Solvability}
%--------------------------------------------------------------------------

The analytical formulas for the Jacobian and the Hessian offer insight into the structure
of the nonlinear RBF system ${\vec W}({\vec\alpha})=0$. Since $\mu$ inherits the smoothness of ${\phi(r)}$, the possible minima of $\mu$ are all critical points.
 
Local minima of the merit function (i.e. those for which $\mu>0$) have the important property that $\det{J}=0$, because $\nabla\mu=J^T{\vec W}=0$ but ${\vec W}\neq 0$. (Obviously, the same property holds for any critical point of $\mu$.) Therefore, nonsingularity of $J({\vec\alpha})$ allows one to guarantee that RBFTrust %with the full TRS method 
will converge to a root--unless it can drift towards $||{\vec\alpha}||\shortrightarrow\infty$ seeking to reduce $\mu$. This is the idea of Theorem \ref{Th:S&U}, for which Theorem \ref{Th:Courant} below will be needed (see \cite{Ghaderi2011} and references therein).

\begin{theorem}[Courant's finite-dimensional mountain pass theorem]

Suppose that $\mu({\vec\alpha})$ is 
\begin{itemize}
\item continuous with continuous derivatives in ${\mathbb R}^N$,
\item coercive (i.e. $\lim\limits_{||{\vec\alpha}||\shortrightarrow\infty} \mu({\vec\alpha})=\infty$), and
\item possesses two different strict (i.e. isolated) minima ${\vec\alpha}_1$ and ${\vec\alpha}_2$.
\end{itemize}
Then, $\mu$ possesses a third critical point ${\vec\alpha_3}$ such that $\mu({\vec\alpha}_1)<\mu({\vec\alpha}_3)>\mu({\vec\alpha}_2)$--therefore distinct from  ${\vec\alpha}_1$ and ${\vec\alpha}_2$.
\label{Th:Courant}
\end{theorem}

\begin{theorem}[Sufficient conditions for solvability and uniqueness of nonlinear RBF collocation]
Let $\mu({\vec\alpha})={\vec W}^T({\vec\alpha}){\vec W}({\vec\alpha})/2$ with ${\vec W}:{\mathbb R}^N\mapsto{\mathbb R}$, and let $J({\vec\alpha})$ be the Jacobian of ${\vec W}$. Assume that $\mu$ is coercive and that $\det{J}\neq 0$ in ${\mathbb R}^N$. Then,
\begin{itemize}
\item There is one unique root ${\vec\alpha}_*$ to ${\vec W}({\vec\alpha})=0$. 
\item Under the further conditions in Theorem \ref{Th:Convergence}, and assuming numerical precision high enough to overcome possible ill-conditioning, RBFTrust with the full TRS method will find it from any initial guess. 
\end{itemize}
\label{Th:S&U}
\end{theorem}

{\em Proof.} If the function is coercive, there is a minimizer of $\mu$ and since $\mu$ is twice differentiable and $J$ is nonsingular, that minimizer is a root. For uniqueness, assume that there were two roots. Then, Theorem \ref{Th:Courant} ensures the existence of a third critical point with $\mu>0$ and therefore singular $J$, which is a contradiction. For the second part, simply apply Theorem \ref{Th:Convergence}. $\square$ \newline

Theorem \ref{Th:S&U} establishes a parallelism between the linear and nonlinear versions of Kansa's method, with $J$ playing the same role as the collocation matrix $K=[{\cal L}\phi]$ from (\ref{F:Kansa}). Unfortunately, the assumptions of Theorem \ref{Th:S&U} are quite elusive--even assuming numerical stability. For instance, consider a linear BVP like (\ref{F:Linear_PDE}), which can be regarded as a particular nonlinear one, and let us analyze its root structure from the point of view of $\mu({\vec\alpha})$. The Jacobian of (\ref{F:Kansa}) is $K$ itself, and $H=K^TK$. If $\det{K}\neq 0$, $H-\big(\lambda_{min}(K)\big)^2I>0$, so that $\mu({\vec\alpha})$ is strongly convex, and thus has a unique minimum, which must be a root. Nonetheless, this could not be proved by Theorem \ref{Th:S&U} since $\lim_{||{\vec\alpha}||\shortrightarrow\infty} \mu({\vec\alpha})=0$ along the direction ${\vec\alpha}=K^{-1}[f({\bf x}_1),\ldots,g({\bf x}_N)]^T$ --the rhs vector in (\ref{F:Kansa})--so that $\mu$ is not coercive in the first place.    Let us now come back to nonlinear problems. Even in the simple and well-behaved Example I in Section \ref{S:Examples}, it does not look trivial to prove (or disprove) coercivity and nonsingularity of $J$ from (\ref{F:J_semilinear}). Summing up, while combining the formulas for $J$ and $H$ with Theorem \ref{Th:S&U} looks theoretically exploitable, we have been unable to do so.

%--------------------------------------------------------------------------
\section{Numerical experiments}
\label{S:Examples}
%--------------------------------------------------------------------------
In this section we test RBFTrust on four nonlinear elliptic BVPs with increasing
level of difficulty, taken from the literature. For comparison purposes, the
operator-Newton method and a Matlab canned implementation are also sometimes
used. In sum, each of the problems I-IV is solved with one or more of
the following methods (see also Table \ref{T:Repex}):
\begin{itemize}
\item The operator-Newton method without smoothing (Section \ref{S:Newton}). We just report whether convergence was or not achieved, mostly with the goal of illustrating the shortcomings of this approach.
\item Matlab's \verb|fsolve| using the TRA with dogleg TRS scheme. This method
relies on a finite-difference approximation to $J$ and, as far as we know, is
the only way that the TRA has been employed in the RBF literature so far.
\item RBFTrust, meaning that one of the three TRS schemes reviewed in Section \ref{S:Trust}
is combined with the analytic approximations to $J$ and $H$ in Section \ref{S:RBF_trust}.
Recall that those three TRS schemes are: dogleg (Section \ref{SS:Dogleg}), full (end
of preamble of Section \ref{S:Trust}), and 2Dsub (Section \ref{SS:2DSub}). In order to compare the various
possibilities, they have been incorporated into a single Matlab code written
from scratch 
%(i.e. not relying on specialized Matlab functions such as \verb|fsolve|)
. Moreover, we sometimes test scaling as described in Section \ref{SS:Scaling}.
RBFTrust variations like those are the new approach introduced and advocated in this paper.
\end{itemize}

With RBFTrust, the RBF formulas for the Jacobian (\ref{F:RBF_Jacobian}) and Hessian (\ref{F:RBF_Hessian}) were
checked against finite-difference approximations produced with DERIVEST (free Matlab program by John D'Errico), and found to agree within the numerical tolerances. The RBFs used in the experiments are those in Table \ref{T:RBFs}.

For every domain $\Omega$, we consider an evaluation set of $N_e>>1$ points (different from collocation nodes) scattered over $\Omega$. The error $\epsilonup$ and the interpolation residual to the PDE, $R$, are monitored there. For instance, the accuracy is estimated via $RMS(\epsilonup):=\sqrt{\sum_{i=1}^{N_e}\epsilonup_i^2/N_e}$. The {\em collocation} residual $R_c$, of course, is evaluated on the collocation nodes.  

Convergence of RBFTrust to a root is declared as soon as $\mu$ stagnates as $\mu=\mu_{\infty}$ and $\mu_{\infty}$ is negligible small. This
is a conservative criterion, for when it happens, $RMS(R)$ and $RMS(R_c)$ will typically
have stagnated some iterations before, and $RMS(\epsilon)$, even before that.
\begin{table}[h!]
\begin{footnotesize}
\caption{Methods used in Problems I-IV (U/S= unscaled/scaled $J$ and $H$).}
\label{T:Repex}
\[\begin{array}{|ll|cccc|}
\multicolumn{2}{c}{} & \textrm{I} & \textrm{II} & \textrm{III} & \multicolumn{1}{c}{\textrm{IV}} \\
\multicolumn{2}{c}{{\bf Method}} & \textrm{{\bf Cubic}} & \textrm{{\bf Plateau's}} &
 \textrm{{\bf Hele-Shaw}} & \multicolumn{1}{c}{\textrm{{\bf Monge-Amp\`ere}}} \\
\hline\noalign{}
\verb|fsolve| & \textrm{dogleg + finite-differences $J$} & \textrm{U} & & & \\
\hline\noalign{}
                  & \textrm{dogleg + analytic $J$}        & \textrm{U} & \textrm{U} & \textrm{U} & \textrm{U} \\
\textrm{RBFTrust} & \textrm{2Dsub + analytic $J$ and $H$} & \textrm{U} & \textrm{U,S} && \textrm{U,S}\\
						& \textrm{full + analytic $J$ and $H$}  & \textrm{U} & \textrm{U,S} && \textrm{U,S}\\
\hline\noalign{}
\textrm{linearization} & \textrm{operator-Newton} &\textrm{U}&\textrm{U}&\textrm{U}&\textrm{U}\\
\hline\noalign{}
\end{array}\]
\end{footnotesize}
\end{table}

\subsection{Example I: diffusion equation with a cubic nonlinearity}

\begin{figure}[h]
\centerline{\includegraphics[width=.6\columnwidth, height=.6\columnwidth]{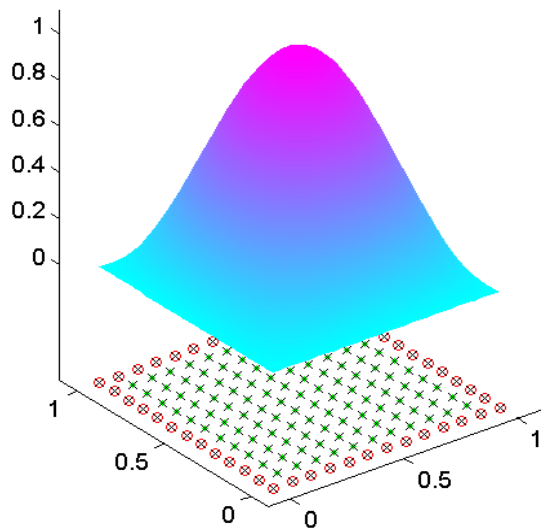}}
\caption{Solution of Example I on an illustrative pointset. Symbols stand for: green $\bullet=$ PDE collocation node; red $o=$ BC collocation node; $\times=$ RBF center.}
\label{I:SolutionI}
\end{figure}

The first example is the semilinear equation (\ref{F:PDE_cubic}) solved on the square $[0,1]^2$ with $f$ such that the exact solution is $u_{ex}(x,y)=\sin{(\pi x)}\sin{(\pi y)}$, and $g=u_{ex}|_{\partial\Omega}=0$ (see Figure \ref{I:SolutionI}). The shrunken $J$ and $H$ are (\ref{F:J_semilinear}) and (\ref{F:H_semilinear}), respectively. We solve this problem on a grid-like pointset. As a starting guess, ${\vec\beta}_0=(0,\ldots,0)$. No scaling was implemented. %Results using two different RBFs--WC4 and MQ--are listed on Table \ref{T:Cubic_with_WC4} and Table \ref{T:Cubic_with_MQ_fsolve}, respectively.   

\begin{table}[h]
\begin{footnotesize}
\caption{$N$ is the number of RBFs, $\mu_\infty$ the value of the merit function at the converged minimum, $\epsilonup$ the error, $\kappa$ the condition number, and the figures in parentheses along the number of iterations are the CPU time in seconds.}\label{T:Cubic_with_WC4}
\[\begin{array}{lcccccc}
\multicolumn{7}{c}{\textrm{Example I solved with RBFTrust and RBF=WC4 ($L=0.3$)}}\\
\hline\noalign{\smallskip}
N 		& \mu_\infty & \textrm{RMS($\epsilonup$)} & \kappa(J) & \textrm{iter(dogleg)} & \textrm{iter(full)}& \textrm{iter(2Dsub)}\\
\hline\noalign{\smallskip}
529	&{\cal O}(10^{-26})	&.01103	&109	&16 (2.18)	&16 (5.35)  &16 (5.34)\\
676	&{\cal O}(10^{-26})	&.00675	&219	&16 (4.49)	&17 (11.31)	&17 (11.20)\\
784	&{\cal O}(10^{-26})	&.00511	&304	&17 (7.11)	&17 (17.60)	&17 (17.50)\\
1024	&{\cal O}(10^{-25})	&.00309	&615	&17 (15.98)	&17 (39.45)	&17 (39.39)\\
1296	&{\cal O}(10^{-25})	&.00204	&1157	&17 (32.50)	&17 (122.72)&17 (73.08)\\
1681	&{\cal O}(10^{-24})	&.00130	&2231	&17 (67.53)	&17 (161.81)&17 (171.64)\\
2116	&{\cal O}(10^{-24})	&.00087	&3982	&18 (144.91)&17 (322.09)&17 (328.72)\\
\hline\noalign{\smallskip}
\end{array}\]
\end{footnotesize}
\end{table}

\begin{table}[h]
\begin{footnotesize}
\caption{Same notation as in Table \ref{T:Cubic_with_WC4}, plus the total number of function evaluations
(fevals). $N = 2116$ not reported because it takes unacceptably long.}\label{T:Cubic_with_WC4_fsolve}
\[\begin{array}{lcccccc}
\multicolumn{7}{c}{\textrm{Example I solved with {\em fsolve} and RBF=WC4 ($L=0.3$)}}\\
\hline\noalign{\smallskip}
N 		& \mu_\infty & \textrm{RMS($\epsilonup$)} & \kappa(J) & \textrm{iterations} & \textrm{CPU time (s.)} & \textrm{fevals}\\
\hline\noalign{\smallskip}
529  &{\cal O}(10^{-25}) &.01103 &109  &4 &308   &2210 \\
676  &{\cal O}(10^{-25}) &.00675 &219  &4 &1063  &2885 \\
784  &{\cal O}(10^{-24}) &.00511 &304  &4 &1446  &3385 \\
1024 &{\cal O}(10^{-24}) &.00309 &615  &4 &4071  &4505 \\
1296 &{\cal O}(10^{-23}) &.00204 &1157 &4 &10577 &5785 \\
1681 &{\cal O}(10^{-23}) &.00130 &2231 &4 &32061 &7610 \\
\hline\noalign{\smallskip}
\end{array}\]
\end{footnotesize}
\end{table}

Let us start by picking Wendland's $WC4$ RBF. On Table \ref{T:Cubic_with_WC4}, all three TRS approximations for RBFTrust converge in a similar
number of iterations, albeit they do not take the same CPU time. This is so because
each TRS approximation takes a different number of function evaluations (each
requiring to solve a linear system), although the ratio evaluations/iteration is
always ${\cal O}(1)$--see Section \ref{S:Trust}. (In fact, 2Dsub ought to be substantially faster than full. The reasons why this does not happen may be that $N$ is not yet large enough, and/or that full uses Matlab's built-in \verb|trust| while 2Dsub is a non-optimized code.)
 
Given the mild nonlinearity of this problem, the operator-Newton method was able to converge from ${\vec\beta}_0$ (as well as from all the other tested initial guesses). Moreover, $H$
is not required in order to enforce convergence, and hence the simplest TRS scheme, dogleg--which
only needs $J$--is the more efficient one. Recall that the three TRS solvers of RBFTrust
in Table \ref{T:Cubic_with_WC4} all use the analytical $J$ and $H$ from Section \ref{S:Trust}. On the other hand,
on Table \ref{T:Cubic_with_WC4_fsolve}, Example I is solved with \verb|fsolve|, whereby every iteration takes
$N-M$ function evaluations ($M$ boundary nodes) in order to approximate $J$ with
finite differences. Since the Jacobians are well conditioned, \verb|fsolve| is able to find the
root the imperfect $J$ notwithstanding. However, computational times are much
longer than with RBFTrust. (Note that $\mu_{\infty}$ is different from before because the
convergence criterion of \verb|fsolve| is independent of that used in RBFTrust.)

Let us now pick the $MQ$ RBF. On Tables \ref{T:Cubic_with_MQ} and \ref{T:Cubic_with_MQ_fsolve}, only the dogleg TRS scheme was used, both with
RBFTrust (i.e. analytical $J$, Table \ref{T:Cubic_with_MQ}) and \verb|fsolve| (i.e. finite-differences approximation
to $J$, Table \ref{T:Cubic_with_MQ_fsolve}). Comparing both, not only does \verb|fsolve| take much
longer, but the finite-differences $J$ is inadequate to provide convergence already
at mild condition numbers (larger than $\kappa(J)\approx 10^8$). After Table \ref{T:Cubic_with_MQ_fsolve}, we will not
show more results involving finite-differences Jacobians.

\begin{figure}[h]
\centerline{\includegraphics[width=.9\columnwidth, height=.45\columnwidth]{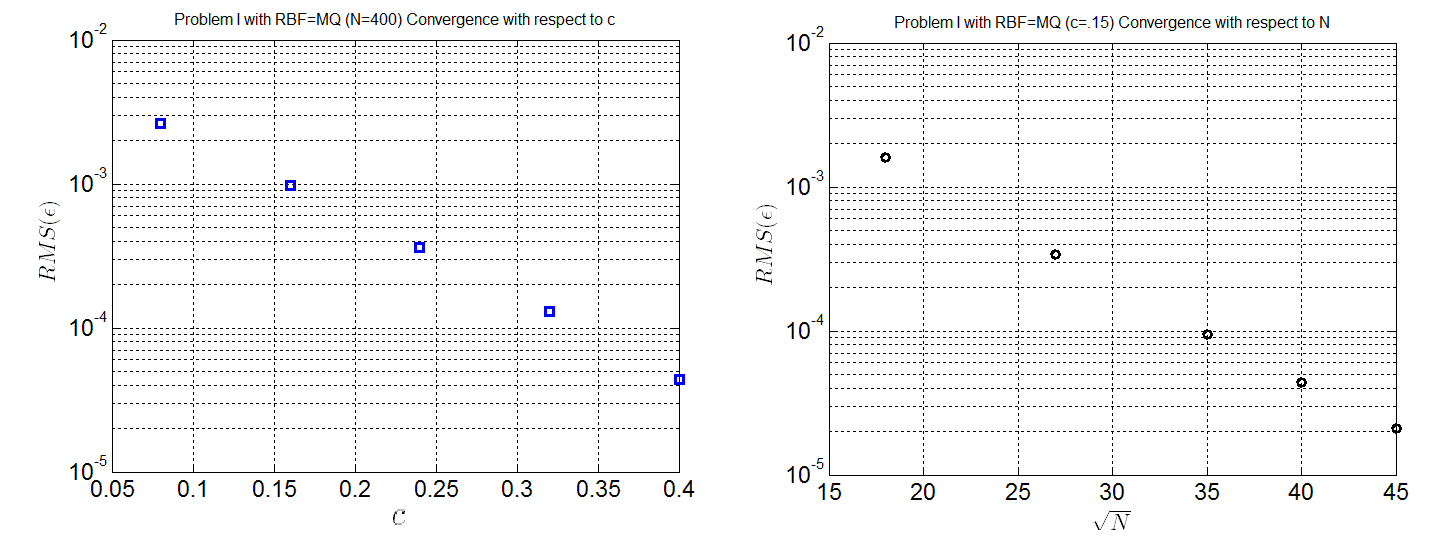}}
\caption{$h-c$ convergence of $RMS(\epsilonup)$ in Example I (data from Table \ref{T:Cubic_with_MQ}). Left: error vs. c. Right: error vs. $\sqrt{N}$ (note that in gridlike pointsets, $1/h\sim\sqrt{N})$}
\label{I:Convergence}
\end{figure}

Also on Table \ref{T:Cubic_with_MQ}, asymptotic exponential convergence--as in linear elliptic
PDEs, see \cite{Cheng2003,Bernal_Newt}--is hinted: $RMS(\epsilonup)= {\cal O}\big( C^{c/h}\big)$, $0<C<1$. The fill distance $h$ is here proportional to $1/\sqrt{N}$ (see Figure \ref{I:Convergence}). 
The possibility of $c-$convergence remains an appealing feature of RBF collocation
also with nonlinear problems.

\begin{table}[h]
\begin{footnotesize}
\caption{Convergence with respect to $c$ and $N$. $R$ is the interpolation residual.}\label{T:Cubic_with_MQ}
\[\begin{array}{lccccccc}
\multicolumn{8}{c}{\textrm{Example I solved with RBFTrust and $RBF=MQ$}}\\
\hline\noalign{\smallskip}
N &  c & \textrm{RMS($\epsilonup$)} & \textrm{RMS($R$)} & \kappa(J) & \mu_\infty & \textrm{iterations} & \textrm{CPU t (s.)}\\
\hline\noalign{\smallskip}
400 &.08 &.00264 &3.65 &{\cal O}(10^2)  &{\cal O}(10^{-25}) &23 &1.07 \\
400 &.16 &.00098 &1.40 &{\cal O}(10^4)  &{\cal O}(10^{-23}) &28 &1.65 \\
400 &.24 &.00036 &0.52 &{\cal O}(10^7)  &{\cal O}(10^{-21}) &28 &1.40 \\
400 &.32 &.00013 &0.19 &{\cal O}(10^9)  &{\cal O}(10^{-18}) &37 &1.71 \\
400 &.40 &4.4\times 10^{-5}&0.07 &{\cal O}(10^{12}) &{\cal O}(10^{-15}) &42 &1.77 \\
\hline\noalign{\smallskip}
324  &.15 &.00160           &  1.86 &{\cal O}(10^3) & {\cal O}(10^{-24})  &26 &0.65 \\
729  &.15 &.00034           &  0.88 &{\cal O}(10^6) & {\cal O}(10^{-21})  &28 &10.3 \\
1225 &.15 &9.5\times 10^{-5}&  0.44 &{\cal O}(10^8) & {\cal O}(10^{-18})  &31 &62.5 \\
1600 &.15 &4.4\times 10^{-5}&  0.28 &{\cal O}(10^{10})& {\cal O}(10^{-17})  &32 &136.7 \\
2025 &.15 &2.1\times 10^{-5}&  0.17 &{\cal O}(10^{11})& {\cal O}(10^{-15})  &33 &296.0 \\
\hline\noalign{\smallskip}
\end{array}\]
\end{footnotesize}
\end{table}

\begin{table}[h]
\begin{footnotesize}
\caption{Same experiments as in Table \ref{T:Cubic_with_MQ}. The case $N = 1225$ was aborted after
$94$ iterations due to exceedingly slow convergence (if at all).}\label{T:Cubic_with_MQ_fsolve}
\[\begin{array}{lcccccccc}
\multicolumn{9}{c}{\textrm{Example I solved with {\em fsolve} and $RBF=MQ$}}\\
\hline\noalign{\smallskip}
N &  c & \textrm{RMS($\epsilonup$)} & \textrm{RMS($R$)} & \kappa(J) & \mu_\infty & \textrm{iterations} & \textrm{CPU time (s.)} & \textrm{fevals}\\
\hline\noalign{\smallskip}
400 &.08 &.00264 &3.65 &{\cal O}(10^2) &{\cal O}(10^{-24}) &6  &117   &2382  \\
400 &.16 &.00098 &1.40 &{\cal O}(10^4) &{\cal O}(10^{-19}) &10 &189   &3575  \\
400 &.24 &.00036 &0.52 &{\cal O}(10^7) &{\cal O}(10^{-14}) &129&1765  &32530  \\
\hline\noalign{\smallskip}
324 &.15 &.00160 &1.86 &{\cal O}(10^3) &{\cal O}(10^{-19}) &8  &80    &2313  \\
729 &.15 &.00034 &0.88 &{\cal O}(10^6) &{\cal O}(10^{-17}) &33 &6956  &20659  \\
1225 &.15& -     & -   &{\cal O}(10^8) &11.98              &94 &131181&78503  \\
\hline\noalign{\smallskip}
\end{array}\]
\end{footnotesize}
\end{table}

\subsection{Example II: Plateau's problem}

\begin{figure}[h]
\centerline{\includegraphics[width=1.2\columnwidth, height=.6\columnwidth]{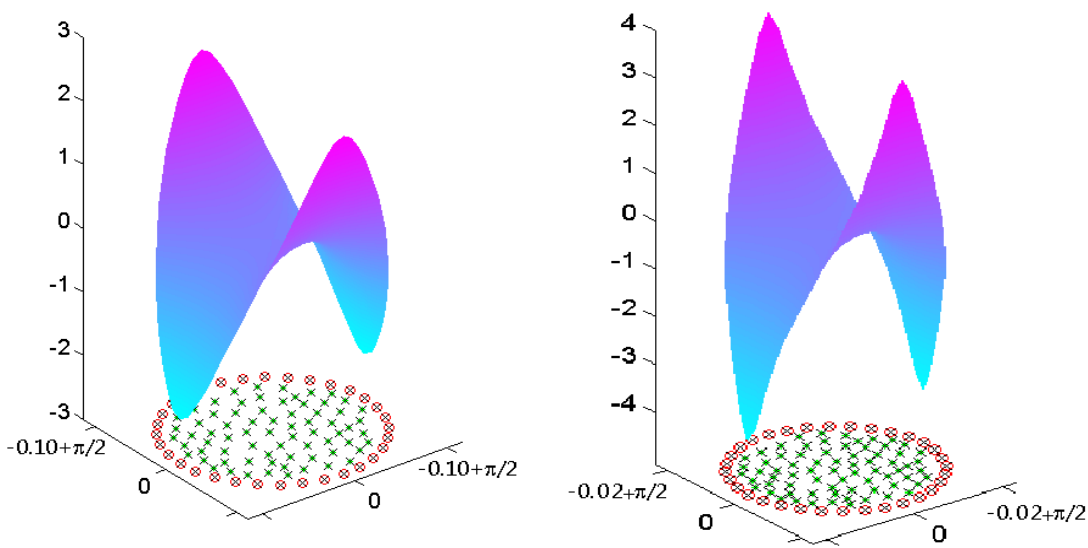}}
\caption{Solutions of Example II with radius $R=-.10+\pi/2$ (left) and $R=-.02+\pi/2$ (right). Note that scales vary.}
\label{I:SolutionsII}
\end{figure}

Plateau's problem (\ref{F:PDE_Plateau}) is a case of the least-surface equation where the solution can
be expressed as a function of $(x,y)$. The solution on the circle 
$\Omega$ centred at the
origin and with radius $R<\pi/2$ is known as Scherk's first minimal surface. Let
$R = \pi/2-s$ with $s>0$. As $s\shortrightarrow 0$, $|\nabla u|$ close to $\partial\Omega$
grows unbounded, thus becoming numerically more challenging--see Figure \ref{I:SolutionsII}.
\begin{eqnarray}
\label{F:PDE_Plateau}
\nabla\cdot\Big( \frac{\nabla u}{\sqrt{1 + |\nabla u|^2}} \Big)= 0 \textrm{ if }{\bf x}\in\Omega,\qquad
u= \log{(\cos{x}/\cos{y})}=u_{ex}|_{\partial\Omega} \textrm{ if }{\bf x}\in\partial\Omega.
\end{eqnarray}

\begin{figure}[h]
\centerline{\includegraphics[width=1.2\columnwidth, height=.6\columnwidth]{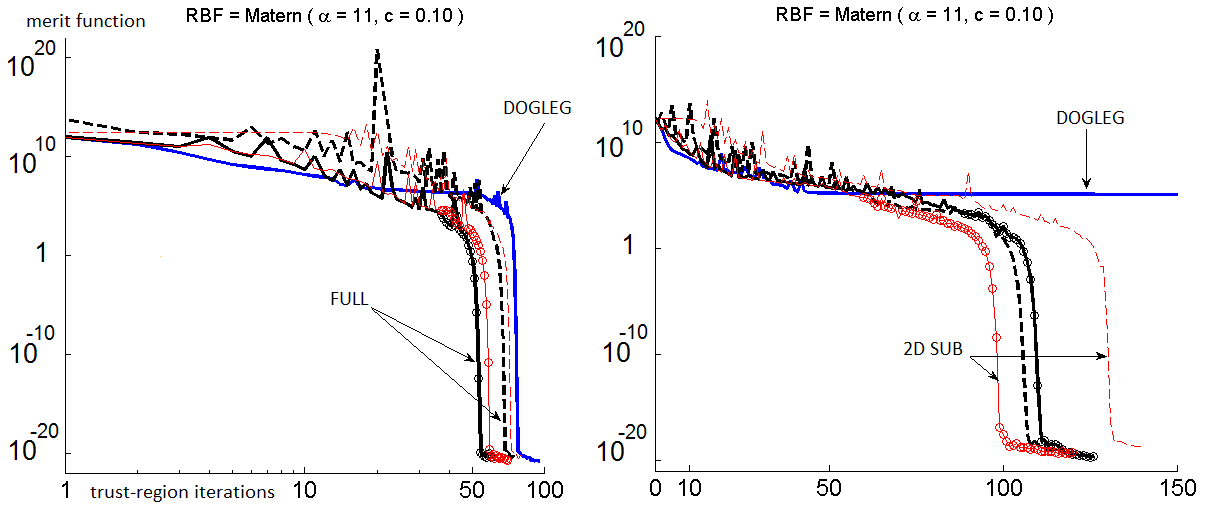}}
\caption{Merit function $\mu$ vs. RBFTrust iterations in Example II (left: radius $R=\pi/2-.10$; right: $R=\pi/2-.02$). TRS methods are: dogleg (blue), full (black) and 2Dsub (red). Dashed lines stand for scaled $J$ and $H$, and circles (only on the unscaled curves) for a positive-definite Hessian at that iteration. The Hessian is critical for finding a root as $s\shortrightarrow 0^+$ (and $|\nabla u|\shortrightarrow\infty$ close to $\partial\Omega)$. Note also that scaled methods tend to take more iterations, which is undesirable.}
\label{I:Plateau}
\end{figure}

We solve two instances of this problem (with $s=.10$ and $s=.02$) with RBFTrust on a scattered pointset with $N=795$ nodes, of which $80$ are along the boundary. The RBF is $MATERN(\alpha=11,c=.10)$. With Example II, we test five TRS methods: dogleg, full (scaled and unscaled), and 2Dsub (scaled and unscaled). The scaling matrix is (\ref{F:Diagonal_scaling}), and $\Delta_0$ is $1$ if unscaled and $10^{10}$ if scaled. The initial guess ${\vec\beta}_0$ is a random vector (but the same for all experiments). Convergence is plotted in Figure \ref{I:Plateau}. In both cases there is apparently a unique root: 
\begin{itemize}
\item [-] For s=.02, ${\vec\alpha}_*:\{RMS(\epsilon)=.00550132,RMS(R)=1021.08,\kappa(J)=112276\}$.
\item [-] For s=.10, ${\vec\alpha}_*:\{RMS(\epsilon)=.00444585,RMS(R)=42.576,\kappa(J)=74966\}$.
\end{itemize}
%${\vec\alpha}_*:\{RMS(\epsilon)=.00550132,RMS(R)=1021.08,\kappa(J)=112276\}$ for $s=0.02$ and ${\vec\alpha}_*:\{RMS(\epsilon)=0.00444585,RMS(R)=42.576,\kappa(J)=74966\}$ for $s=0.10$. 
While for $s=.10$ all five TRS variants of RBFTrust correctly converge to the root, in the more difficult case ($s=.02$) dogleg stagnates at $\mu>>0$ (with $RMS(\epsilon)\approx .443$). Actually, this value is not a local, non-root minimum (in fact $\kappa(J)=275042$), but rather the dogleg steps yield the same very small drop as the Cauchy steps ($\approx 193$). Therefore, Hessian information is absolutely critical in order to solve Example II with $s=.02$. Both the full and 2Dsub use $H$. However, since 2Dsub solves the TRS in a two-dimensional, rather than $(N-M)-$dimensional, space, each 2Dsub iteration is much cheaper than a full iteration. Therefore, 2Dsub solves the problem faster, even though it may take more iterations.  

We also stress the robustness of RBFTrust, spanning $33$ orders of magnitude (from $\mu({\vec\alpha}_0)\approx 10^{12}$ to $\mu({\vec\alpha}_\infty)\approx 10^{-21}$) over a non-convex merit function landscape and from a Gaussian random guess. Many more unreported experiments show that the qualitative results in Figure \ref{I:Plateau} hold for other RBFs and for other random guesses, always finding the same RBF solution. 

On the other hand, scaling proves a disappointment, and in three out of four cases scaled versions take longer to converge. This is seemingly due to the worsening of ${\kappa(H)}$.\newline

\begin{table}[h]
\begin{footnotesize}
\caption{Example II (with $s=1$) solved with the operator-Newton method from a perturbed Laplacian guess. Left: no perturbation. Middle: convergence still takes place. Right: a slightly larger perturbation leads to divergence.}\label{T:opNewt}
\[\begin{array}{|c|ll|ll|ll|}
\multicolumn{1}{c}{} & \multicolumn{2}{c}{\textrm{$\vec{\alpha}_0$=Laplacian}} &
\multicolumn{2}{c}{\textrm{$\vec{\alpha}_0$=Laplacian + $\vec{\delta}$/250}} &
\multicolumn{2}{c}{\textrm{$\vec{\alpha}_0$=Laplacian + $\vec{\delta}$/247}} \\
\hline
\textrm{iteration $\sharp$} & \textrm{RMS($\epsilonup$)} & \textrm{RMS($R_c$)} & \textrm{RMS($\epsilonup$)} & \textrm{RMS($R_c$)} & \textrm{RMS($\epsilonup$)} & \textrm{RMS($R_c$)}\\
\hline
0 & .00124              & .1464	            &	.00124	         & 0.1464	&0.00124& 0.1464 \\
1 & .000378             & .059               &  .0603	            & 2.29	&0.1221 & 3.93 \\
2 & .000162             & .030               &	.0067	            & 1.18	&0.0594 & 5.42 \\
3 & 5.00\times 10^{-5}  & .0151 	            &  .0020	            & 0.4053 &0.0781 & 17.29 \\
4 & 4.20\times 10^{-5}  & .0081	            &  .00071            & 0.2001 &0.5668 & 19185 \\
10 & 1.43\times 10^{-5} & .00024             &	1.73\times 10^{-5}& 0.0041	& \infty&\infty\\ 
60 & 1.39\times 10^{-5} & 6.6\times 10^{-16} &	1.39\times 10^{-5}& 6.15\times 10^{-15}	& \infty&\infty\\
\hline\noalign{\smallskip}
\end{array}\]
\end{footnotesize}
\end{table}

Let us use Example II to illustrate the shortcomings of the operator-Newton method. As expected, and unlike RBFTrust before, the operator-Newton method was unable to find a root from a random guess. Therefore, we picked ${\vec\beta}_0$ from ${\vec\alpha}_0$ resulting of interpolating the solution of a Laplace equation with the same BCs over the same pointset (we call this the "Laplacian guess"). Even with the Laplacian guess, the radius of $\Omega$ must be substantially smaller than before in order to get the operator-Newton method to converge. Consequently, we kept the same RBFs and $N$ but shrank $\Omega$ to $R=\pi/2-1$(i.e. we let $s=1$). Notice that the solution is much flatter now and, in sum, the BVP is much less challenging than the one solved before with RBFTrust. Only after those simplifications, the operator-Newton method was finally able to converge. Further investigating the performance of the operator-Newton, let $\vec{\delta}$ be a random vector drawn from a standard normal distribution, ${\cal N}(0,1)$. Table \ref{T:opNewt} shows the effect on convergence of a small perturbation over the Laplacian guess, highlighting the lack of robustness of the operator-Newton method. Notice that the residual $R_c$ in Table \ref{T:opNewt} is the {\em collocation} residual, and that values of RMS($R_c$)$\lesssim 10^{-14}$ strongly hint to a root of the collocation system (\ref{F:Nonlinear_System}). 

Wrapping up Example II, the operator-Newton method performs very poorly compared to RBFTrust with dogleg, and the latter worse than Hessian-based RBFTrust. In fact, in the presence of high gradients, even dogleg with analytic $J$ fails, and second-order information derived from the analytic $H$ becomes indispensable.

\subsection{Example III: simulation of powder injection molding}

The three-dimensional creeping flow of molten polymer into a thin cavity can in many cases be modeled by a two-dimensional free-boundary problem involving the p-Laplace operator. This is known as the Hele-Shaw approximation \cite{Bernal2007,Bernal09}. At every timestep, the following nonlinear elliptic BVP for the pressure field $u(x,y)$ must be solved (here in dimensionless units): 

\begin{eqnarray}
\label{F:PIM}
\nabla\cdot\Big( |\nabla u|^{\gamma} \nabla u\Big)= 0
\textrm{ if }{\bf x}\in\Omega,\qquad
\textrm{with mixed BCs}
\left\{
\begin{array}{ll}
u= 1 & \textrm{ if } {\bf x}\in\Gamma_I\\
u= 0 & \textrm{ if } {\bf x}\in\Gamma_F\\
\partial u/\partial N=0 & \textrm{ if } {\bf x}\in\Gamma_W.
\end{array}
\right.
\end{eqnarray}

\begin{figure}[h]
\centerline{\includegraphics[width=1.2\columnwidth, height=.6\columnwidth]{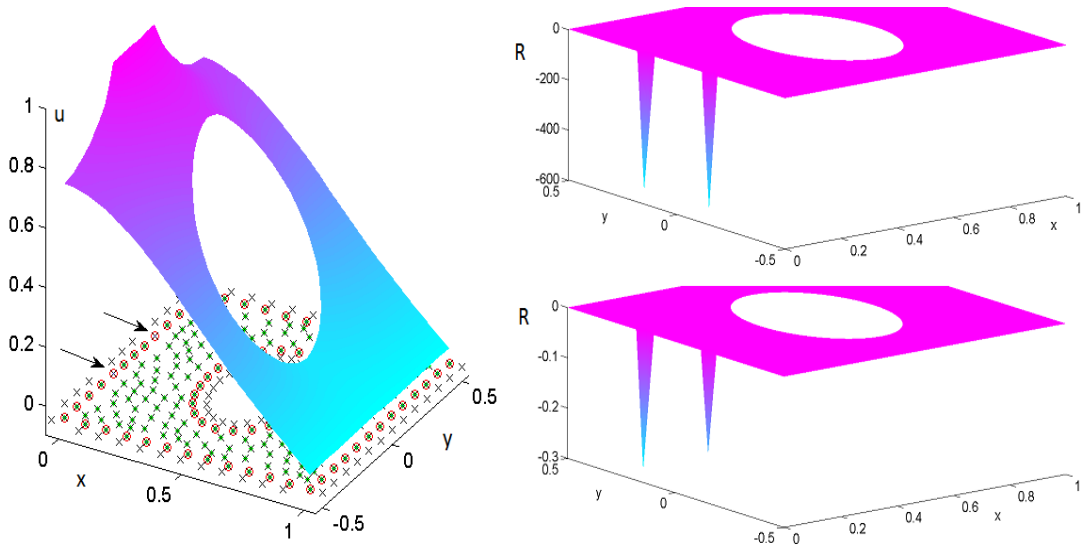}}
\caption{Left: Solution of Example III. The injection inlet $\Gamma_I$ is between the arrows. Notice the extra RBF centres ($\times$) placed outside for PDEBC; except at both ends of $\Gamma_I$--where only the Dirichlet BC is enforced, for $u$ is not smooth there. Along the front $\Gamma_F$, $u=0$. Right: residual peaks at the singularities with $0+0$ Motz functions (up) and $2+2$ (bottom). Note that the vertical scale varies.}
\label{I:PIM}
\end{figure}

$\Gamma_I$ represents the injection slit (where the pressure is enforced by the injection machine), $\Gamma_F$ is the (frozen) free boundary, $\Gamma_W$ are the walls of the floor view of the mold, and $\Gamma_I\cup\Gamma_F\cup\Gamma_W=\partial\Omega$. A typical value is $\gamma=0.6$, which models polyethylene. Figure \ref{I:PIM} (left) shows a FEM numerical solution obtained over a very fine mesh, which we take as a reference.

In addition to the nonlinearity, this problem has two more challenging features. The first one are the BCs of derivative type, which degrade the quality of the RBF solution. The second issue is the singularity in $u$ which takes place on both ends of $\Gamma_I$ due to the change of type of the BCs. Since the RBF interpolant is made up of infinitely smooth functions, it cannot reproduce the singularities and brings about oscillations around them, resulting in large residual peaks--see Figure \ref{I:PIM} (right).
 
We counter the former difficulty by enforcing the PDE also on some (here, all but two) of the boundary nodes. This is the so-called PDEBC strategy \cite{Fedoseyev2002}, whereby as many additional RBFs are added (usually placed off the boundary) to keep the system square--see Figure \ref{I:PIM} (left). Let us now address the presence of nonsmooth boundary singularities. In the case of Laplacian flows (i.e. $\gamma=0$), the singularities can be effectively captured by enriching the RBF interpolant as in (\ref{F:Enriched_RBF_Interpolant}) with Motz functions \cite{Bernal_Newt}:

\begin{equation}
\label{F:Motz}
h_k(r,\theta)= r^{(2k-1)/2}\cos \big[\big( \frac{2k-1}{2}\big)\theta \big],\qquad k\geq 1.
\end{equation} 

For $\gamma=0$, Motz functions (\ref{F:Motz}) do not contribute residual, because they are harmonic. Table \ref{T:PIM} shows the results of testing the same idea with the p-Laplacian flow (\ref{F:PIM}). The collocation pointset is made up of $N=1152$ scattered nodes ($166$ along the boundary). Note that now there are $1152+164$ RBFs due to PDEBC (the PDE is not enforced on the singularities, but the Dirichlet BCs are). Two RBFs were used, the MQ and the IMQ. The initial guess is the solution of a Laplace equation with the same BCs. While the overall effect is not as dramatic as in the linear case, adding a few ($n$) Motz functions still shows its benefits as the condition numbers grow ($n=a+a$ means that $a$ Motz functions (\ref{F:Motz}) with $k=1,\dots,a$ are centred on each singularity). Without enrichment, accuracy does not improve out of a larger value of $c$ (it actually worsens). Moreover, iterations are fewer, and the derivatives of the solution close to the injection area are also much improved, which is important if, for instance, the velocity field is needed. On the other side, because Motz functions {\em do} contribute residual when $\gamma\neq 0$, there seems to be an optimal number of them, depending on the RBF discretization.

\begin{table}[h]
\caption{Effect of $n$ Motz functions on the RBF interpolant in Example III.}\label{T:PIM}
\begin{footnotesize}
\[\begin{array}{lccccccc}
\multicolumn{8}{c}{\textrm{RBF=IMQ ($c=.75$). Solver: RBFTrust+dogleg TRS}}\\
\hline\noalign{\smallskip}
n & \textrm{iter} & \mu_\infty & \textrm{RMS($\epsilonup$)} & \textrm{MAX($\epsilonup$)} & \textrm{RMS($R$)} & \textrm{MAX($R$)} & \kappa(J)\\
\hline\noalign{\smallskip}
0+0&	16	&1.06\times 10^{-22}	&.0028	&.015	&20.35&497.07	&1.59\times 10^{10} \\
1+1&	15	&2.01\times 10^{-23}	&.0101	&.020	&1.94	&59.11	&4.25\times 10^{10} \\
3+3&	14	&5.84\times 10^{-24}	&.0032	&.009	&0.23	&6.14		&8.37\times 10^{10} \\
6+6&	16	&3.60\times 10^{-23}	&.0027	&.014	&0.06	&1.94		&3.88\times 10^{10}\\
9+9&	19	&5.77\times 10^{-22}	&.0034	&.013	&1.06	&25.96	&1.63\times 10^{10} \\
\hline\noalign{\smallskip}
\multicolumn{8}{c}{\textrm{RBF=IMQ ($c=1.50$). Solver: RBFTrust+dogleg TRS}}\\
\hline\noalign{\smallskip}
n & \textrm{iter} & \mu_\infty & \textrm{RMS($\epsilonup$)} & \textrm{MAX($\epsilonup$)} & \textrm{RMS($R$)} & \textrm{MAX($R$)} & \kappa(J)\\
\hline\noalign{\smallskip}
0+0	&49	&1.35\times 10^{-16}	&.0051	&.016	&22.69&547.59&1.32\times 10^{11} \\
1+1	&25	&9.90\times 10^{-19}	&.0018	&.005	&0.01 &0.27	&2.37\times 10^{11} \\
2+2	&20	&9.28\times 10^{-19}	&.0018	&.005	&0.01	&0.28 &2.73\times 10^{11}\\
3+3	&23	&1.47\times 10^{-18}	&.0028	&.005	&0.01	&0.40	&1.02\times 10^{13}\\
4+4	&22	&9.66\times 10^{-19}	&.0033	&.008	&0.02	&0.55	&1.31\times 10^{14}\\
\hline
\end{array}\]
\end{footnotesize}
\end{table}

Despite the remarkable "filing" of the residual peaks by the Motz functions, they do not entirely vanish (see Figure \ref{I:PIM}, bottom right). For that reason, the operator-Newton method also fails to solve this problem, even if Motz functions are incorporated into the interpolant (results not shown). 

Example III was solved exclusively with the RBFTrust+dogleg method. The reason is that condition numbers are higher now and they actually worsen with $n$ (see Table \ref{T:PIM}). The dogleg TRS scheme is the only one which can handle $A$ matrices with such condition numbers, thanks to the fact that for it, $A=J^TJ$--see item 3 in Section \ref{S:Discussion} for further details.

\subsection{Example IV: Monge-Amp\`ere equation in 3D}

The numerical handling of Monge-Amp\`ere (and, in general, of second-order fully nonlinear equations) is considered very difficult. According to \cite{Feng09}, the following was the only numerical example of a Monge-Amp\`ere PDE in 3D at the time of publication:

\begin{eqnarray}
\label{F:MongeAmpere_in_3D}
\det(D^2_3 u)= f = (1+r^2)\exp{(3r^2/2)} \textrm{ if }{\bf x}\in\Omega,\qquad
u= \exp{(r^2/2)} \textrm{ if }{\bf x}\in\partial\Omega,
\end{eqnarray}

where $r=\sqrt{x^2+y^2+z^2}$. The domain is the unit cube $[0,1]^3$ with $u_{ex}$ being the same as the BC. When $f>0$ (as here), (\ref{F:MongeAmpere_in_3D}) has a unique convex solution (for which the PDE is elliptic), but may have other non-convex solutions \cite{Cheng1977}. %For $f=(1+r^2)\exp{(3r^2/2)}$, the unique convex solution is $u_{ex}(x,y,z)=\exp{(r^2/2)}$, where $r=\sqrt{x^2+y^2+z^2}$. 

\begin{figure}[h]
\centerline{\includegraphics[width=1.2\columnwidth, height=.6\columnwidth]{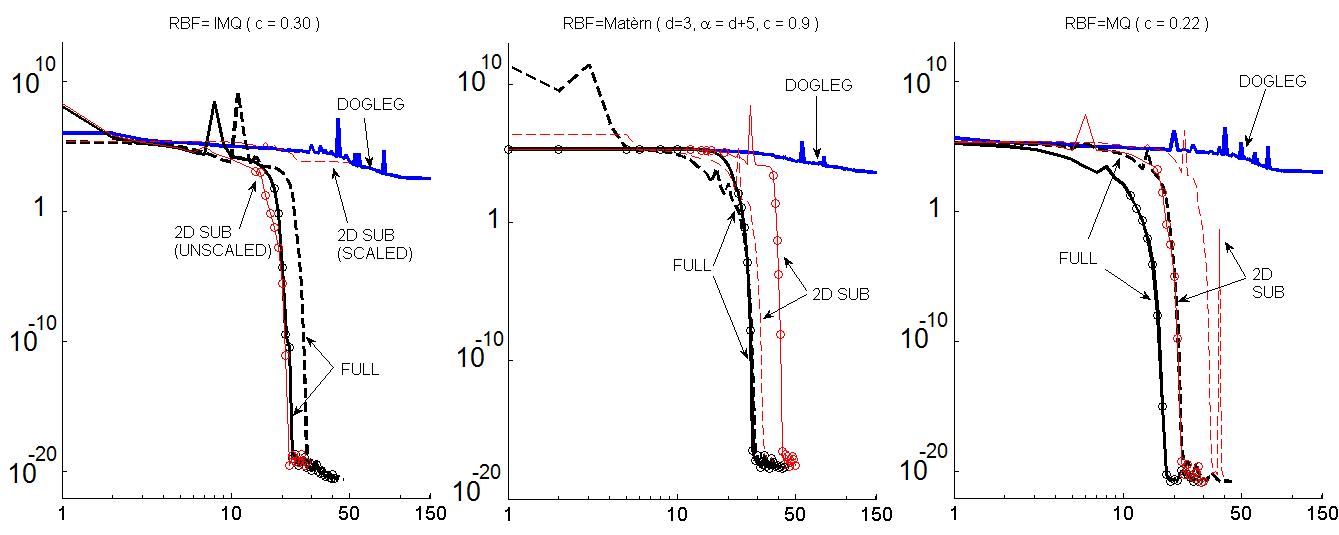}}
\caption{Merit function vs. RBFTrust iterations in Example IV, with $N=2197$ and using $3$ RBFs: IMQ$(c=.30)$ (left), MATERN$(\alpha=8,c=.90)$ (centre), MQ($c=.22$) (right). Symbols as in Figure \ref{I:Plateau}. Again, $H$ is vital for convergence.}
\label{I:MongeAmpere}
\end{figure}

We solved (\ref{F:MongeAmpere_in_3D}) on a lattice-like pointset with $N=2197$ nodes (of which $866$ are on $\partial\Omega$), using the solution of the related Poisson equation as a starting guess. The highly nonlinear character of this PDE is reflected in the fact that, for $N\gtrsim 900$ RBFs, RBFTrust with dogleg TRS ceases to converge to a root. For larger problems than that, the Hessian is needed in order to obtain convergence. This pattern is confirmed by all three RBFs used on Figure \ref{I:MongeAmpere}.

In Figure \ref{I:MongeAmpere}, the dogleg iterations end up having a similar (and extremely small) convergence rate with steepest descent and, for practical purposes, get all but stuck. On the other hand, using $H$ results in $\mu$ dropping by several orders of magnitude the value of the Cauchy step, and eventually in convergence. For each of the RBFs in Figure \ref{I:MongeAmpere}, the full and 2Dsub TRS schemes converge to the same minimum, suggesting that RBFTrust is correctly picking the unique convex solution (see Table \ref{T:MongeAmpere}).

In Example IV, scaling of 2Dsub and full (dashed lines on Figure \ref{I:MongeAmpere}; implemented as in Example II) had again a mixed performance, taking sometimes more iterations than the unscaled versions of RBFTrust and even failing to converge with the MQ RBF (red dashed curve on the leftmost part on Figure \ref{I:MongeAmpere}). Seemingly, diagonal scaling according to (\ref{F:Diagonal_scaling}) was destabilizing RBFTrust.% This serves as a warning that the scaling recipe (\ref{F:Diagonal_scaling}) is not foolproof. 

As was to be expected, the operator-Newton method systematically failed to solve Example IV except for small values of $N$.

\begin{table}[h]
\caption{Properties of ${\vec\alpha}_{\infty}$ in Problem IV for several RBFs.}\label{T:MongeAmpere}
\begin{footnotesize}
\[\begin{array}{lccccc}
\hline\noalign{\smallskip}
RBF &  \mu_{\infty} & \textrm{RMS($\epsilonup$)} & \textrm{RMS($R$)} & \kappa(J) & \kappa(H)\\
\noalign{\smallskip}\hline\noalign{\smallskip}
\textrm{Mat\'ern} (\alpha=12,c=.80) & 3.2\times 10^{-20}& .00760& 121.51& 202481 & 4.10\times 10^{10} \\
\textrm{IMQ} (c=.30) & 2.5\times 10^{-21} &	.00187	&	26.61 & 1.16\times 10^{6} & 1.36\times 10^{12}\\
\textrm{MQ} (c=.22) & 1.8\times 10^{-21}	&.00027	&	6.03 & 173323 & 3.00\times 10^{10}\\
\textrm{Mat\'ern} (\alpha=8,c=.90) & 2.6\times 10^{-18}&5.08\times 10^{-5}&1.13 &14828 &2.20\times 10^{8}\\
\hline
\end{array}\]
\end{footnotesize}
\end{table}

Table \ref{T:MongeAmpere} shows that the final accuracy (and conditioning) of the root is quite dependent on the RBF used, and roughly proportional to the final interpolation residual. We underscore that all accuracies on Table \ref{T:MongeAmpere} are at least one order of magnitude better that those reported in \cite{Feng09}.

\section{Discussion of the results of problems I-IV}\label{S:Discussion}

%Let us summarize the general patterns arising from Examples I-IV:

	\begin{itemize}
	\item[1.] In all well-conditioned problems which we have run, we have found a distinct minimum of the merit function $\mu$ below or around the machine error ($\approx 10^{-15}$). In such cases, we observe that, at the converged minimum: i) $\kappa(J)$ clearly is not numerically singular; ii) $H>0$; and iii) the final rate of convergence is quadratic. For those reasons we tend to believe that those converged minima are, in fact, roots of the collocation system. Moreover, those roots are independent of the initial guess or TRS method used, strongly suggesting that they are, in fact, unique. %The same observations apply in most cases at the border of numerical instability. 

	\item[2.]Not in a single numerically stable experiment, whether reported or not, did we find evidence of multiple roots. 

	\item[3.] The tradeoff between accuracy and stability--a hallmark of Kansa's method--carries over to nonlinear equations, albeit via a different mechanism. Picking RBF interpolation spaces with better approximation properties (by letting $N$ or the RBF shape parameter grow), pushes the condition number of the matrices involved ($J$ and $A$) towards numerical singularity. This, in turn, slows down the convergence rate of RBFTrust and lifts the value of the minimum of the merit function that can be resolved. Beyond the threshold of numerical instability, the RBFTrust iterations get stuck in a local minimum with much residual and poor accuracy.

This tradeoff has a bearing on the better choice of the TRS approximation method when higher accuracy is desired. On the one hand, $A=H$ (in the full and 2Dsub methods) should entail a better modeling of the complicated residual landscape and thus allow deeper steps to be taken, resulting in fewer iterations and better probability to find the root. On the other hand, $A=J^TJ$ (in the dogleg method) is numerically advantageous because inverting $J^TJ{\vec x}={\vec b}$ can be split into $J^T{\vec y}={\vec b}$ and then $J{\vec x}={\vec y}$. Assuming--from (\ref{F:Gradient_and_Hessian})--that $\kappa(H)\approx\kappa(J^TJ)=\kappa^2(J)$, the dogleg method can in effect yield better RBF approximations of the solution before becoming unstable than would be possible using the full Hessian.

	\item[4.] The performance of the RBFTrust with the 2Dsub TRS solver is also remarkably affected by ill-conditioning. We observed that for Hessians with $\kappa(H)\gtrsim 10^{10}$, the estimation of the smallest algebraic eigenpair via Matlab {\em eigs} failed because the Lanczos iterations did not converge. In such cases, in order to test the 2Dsub scheme, we computed the exact eigenpair via the QR decomposition of $H$--thus defeating the purpose of speeding up the computations. In order to fully exploit the two-dimensional TRS approximation within RBFTrust, the Lanczos method must be replaced by another algorithm more resistent to RBF ill conditioning.

	\item[5.] More often than not, diagonal scaling (\ref{F:Diagonal_scaling}) proved to be a counter--productive addition which compounded the conditioning of $H$.

	\end{itemize}

%--------------------------------------------------------------------------
\section{Conclusions}
\label{S:Conclusions}
%--------------------------------------------------------------------------
The RBFTrust approach developed here is a robust, easy to code, fast and accurate solver for quite general nonlinear elliptic equations. Experiments show that it is vastly superior to previous approaches based on the operator-Newton method or the dogleg method with finite-difference Jacobians. Particularly critical is the use of the analytic formulas for the Jacobian and Hessian of the collocation system made available in this paper--especially so when the equation is highly nonlinear. 

%The two-dimensional subspace approximation of Byrd, Schnabel,... has been tested with success. In order to fully exploit it, a Lanczos-like algorithm tailored to typical RBF matrices--full and ill-conditioned--is still needed. %As usual in meshless collocation methods, accuracy is essentially limited by stability. 

This paper deals exclusively with strict (Kansa-like) RBF collocation. %We have shown that the question of the existence of a strict root of the collocation system in linear and nonlinear BVPs are connected. 
While the numerical results are very good, the theoretical considerations about solvability and uniqueness in Section \ref{S:Solvability} were left undecided in either sense. In practice, anyways, strict collocation is moot as long as the converged residuals are negligible. %This opens up the application of the RBF/TRA approach to more general RBF formulations based on least-squares rather than collocation. 
The ultimate goal is to tackle much larger problems, along the lines of the RBF-PU method \cite{Larsson2015}.% In such cases, conjugate-gradient methods for the TRS warrant closer attention.    

%\end{flushleft}

\section{Acknowledgements}
Portuguese FCT funding under grant SFRH/BPD/79986/2011 and a KAUST Visiting Scholarship at OCCAM in %the Department of Mathematics of 
the University of Oxford are acknowledged.%\newline

The author thanks Holger Wendland for suggesting Example IV and for helfpul discussions while in Oxford. The referees are acknowledged for their meticulous reading, which led to a better and clearer presentation of the paper.

%----------------------------------------------------------------------------------
% BIBLIOGRAPHY
%----------------------------------------------------------------------------------

\end{document}